\documentclass[a4paper,11pt]{article}

{}{}\usepackage[dvips]{graphics}
{}{}\usepackage[matrix,arrow,tips,curve]{xy}
\usepackage[english]{babel}
\usepackage[leqno]{amsmath}
\usepackage{amssymb,amsthm}
\usepackage{enumerate}
{}{}\usepackage{psfrag}
{}{}\psfrag{GC}{\fontsize{20}{20}$\Gamma_C$}
{}{}\psfrag{1}{\fontsize{20}{20}$1$}
{}{}\psfrag{2}{\fontsize{20}{20}$2$}
{}{}\psfrag{3}{\fontsize{20}{20}$3$}
{}{}\psfrag{C1}{\fontsize{20}{20}$C_1$}
{}{}\psfrag{C2}{\fontsize{20}{20}$C_2$}
{}{}\psfrag{boh}{}
{}{}\psfrag{n1}{\fontsize{20}{20}$n_1$}
{}{}\psfrag{n2}{\fontsize{20}{20}$n_2$}
{}{}\psfrag{n3}{\fontsize{20}{20}$n_3$}
{}{}\psfrag{p1}{\fontsize{20}{20}$p_1$}
{}{}\psfrag{p2}{\fontsize{20}{20}$p_2$}
{}{}\psfrag{p3}{\fontsize{20}{20}$p_3$}
{}{}\psfrag{q1}{\fontsize{20}{20}$q_1$}
{}{}\psfrag{q2}{\fontsize{20}{20}$q_2$}
{}{}\psfrag{q3}{\fontsize{20}{20}$q_3$}
{}{}\psfrag{E}{\fontsize{20}{20}$E$}
{}{}\psfrag{U}{\fontsize{20}{20}$U$}
{}{}\psfrag{n}{\fontsize{20}{20}$n$}
{}{}\psfrag{C}{\fontsize{20}{20}$C$}
{}{}\psfrag{X}{\fontsize{20}{20}$X$}
{}{}\psfrag{wX}{\fontsize{20}{20}$\widetilde{X}$}
{}{}\psfrag{E1}{\fontsize{20}{20}$E_1$}
{}{}\psfrag{E2}{\fontsize{20}{20}$E_2$}
{}{}\psfrag{E3}{\fontsize{20}{20}$E_3$}
{}{}\psfrag{Ej}{\fontsize{20}{20}$E_j$}
{}{}\psfrag{Ek}{\fontsize{20}{20}$E_k$}
{}{}\psfrag{pj}{\fontsize{20}{20}$p_j$}
{}{}\psfrag{pk}{\fontsize{20}{20}$p_k$}
{}{}\psfrag{qj}{\fontsize{20}{20}$q_j$}
{}{}\psfrag{qk}{\fontsize{20}{20}$q_k$}
{}{}\psfrag{p}{\fontsize{20}{20}$p$}
{}{}\psfrag{q}{\fontsize{20}{20}$q$}
{}{}\psfrag{ni}{\fontsize{20}{20}$n_i$}
{}{}\psfrag{XII}{\fontsize{20}{20}$X_{j,k}$}
{}{}\psfrag{XIII}{\fontsize{20}{20}$X_{1,2,3}$}
{}{}\psfrag{a}{\fontsize{20}{20}$\alpha$}
{}{}\psfrag{b}{\fontsize{20}{20}$\beta$}
{}{}\psfrag{c}{\fontsize{20}{20}$\gamma$}
{}{}\psfrag{pr2}{\fontsize{20}{20}$h^{-1}(L_1,L_2)\simeq\mathbb{P}^2$}

\newtheorem{thm}{Theorem}[subsection]
\newtheorem{lemma}[thm]{Lemma}
\newtheorem{prop}[thm]{Proposition}
\newtheorem{cor}[thm]{Corollary}
\newtheorem{lemmadef}[thm]{Lemma - Definition}

\theoremstyle{definition}
\newtheorem{defi}[thm]{Definition}

\theoremstyle{remark}
\newtheorem{remark}[thm]{Remark}
\newtheorem{example}[thm]{Example}

\newcommand{\la}{\longrightarrow}
\newcommand{\w}{\widetilde}
\newcommand{\ov}{\overline}
\newcommand{\un}{\underline}
\newcommand{\ma}{\mathcal}
\def\ev{\Delta}
\def\X{\mathcal X}

\def\N{\mathcal N}
\newcommand{\g}{\Gamma}
\newcommand{\ph}{\varphi}
\newcommand{\pr}[1]{\mathbb{P}^{#1}}
\newcommand{\ol}{\mathcal{O}}
\newcommand{\A}{\mathbb{A}}
\newcommand{\Z}{\mathbb{Z}}

\newcommand{\C}{\mathbb{C}}
\newcommand{\Xn}{X^{\nu}}

\newcommand{\Aut}{\operatorname{Aut}}
\newcommand{\Supp}{\operatorname{Supp}}

\newcommand{\Pic}{\operatorname{Pic}}
\newcommand{\Hilb}{\operatorname{Hilb}^{d,g}_{\pr{s}}}

\newcommand{\Id}{\operatorname{Id}}
\newcommand{\im}{\operatorname{Im}}
\newcommand{\Tw}{\operatorname{Tw}}

\newcommand{\mg}{M_g}
\newcommand{\mgbar}{\ov{M}_g}

\newcommand{\mgnbar}{\ov{M}_{g,n}}
\newcommand{\snc}{\ov{S}^{\,r}_C(N)}
\newcommand{\sn}{{S}^{\,r}_f(\mathcal{N})}
\newcommand{\snbar}{\ov{S}^{\,r}_f(\mathcal{N})}
\newcommand{\snf}{\ov{\mathcal{S}}^{\,r}_f(\mathcal{N})}
\newcommand{\sgnf}{\ov{\mathcal{S}}^{\,r,\,l,\,\un{m}}_{g,n}}
\newcommand{\sgnbar}{\ov{S}^{\,r,\,l,\,\un{m}}_{g,n}}
\newcommand{\sgn}{{S}^{\,r,\,l,\,\un{m}}_{g,n}}
\newcommand{\sgbar}{\ov{S}^{\,r,\,l}_{g}}
\newcommand{\sg}{{S}^{\,r,\,l}_{g}}
\newcommand{\sgtil}{\widehat{S}^{\,r,\,l}_g}
\newcommand{\sC}{\widehat{S}^{\,r,\,l}_C}
\def\pdg{\overline{P}_{d,\,g}}

\newcommand{\fc}{\widehat{S}^{\,3,\,0}_C}
\newcommand{\dnmodr}{[\,\un{d}\,]_{r}}
\newcommand{\PC}{\overline{P}_C^{\,d}}
\newcommand{\PCd}{\overline{P}_C^{\ \un{d}}}
\newcommand{\jar}{{R}_g^{1/r}(\omega^{\otimes l})}
\newcommand{\sbarC}{\ov{S}^{\,r,\,l}_C}

\begin{document}
\begin{center}
{\bf\large
MODULI OF ROOTS OF LINE BUNDLES ON CURVES} \\

\bigskip

Lucia Caporaso, Cinzia Casagrande, Maurizio Cornalba

\bigskip

\end{center}

\bigskip

{\small
\noindent
{\bf Abstract:}
We treat the problem of completing the moduli space for roots of line
bundles on curves.  
Special attention is devoted to higher spin curves within the
universal Picard scheme. Two new 
different constructions, both using line bundles on nodal curves as boundary points, 
are carried out and compared with pre-existing ones.

\

\

}

{\small \tableofcontents}
\section{Introduction}
The problem of compactifying the moduli space for roots of line 
bundles on
curves is the object of
this paper.

More precisely, consider a family of nodal curves $f\colon\ma{C} \to 
B$,
with nonsingular fibers over an open subset $U$ of $B$; let $\N$ be a 
line bundle on $\ma{C}$,
viewed as a family of line bundles on the fibers of $f$.
Assume that the relative degree of $\N$ is divisible by some positive 
integer $r$.
Then the family of $r$-th roots\footnote{If $L$ and $N$ are line 
bundles on
$X$ such that $L^{\otimes r}=N$ we call $L$ an $r$-{\it th root} of 
$N$.}
of the restriction of
$\N$ to the nonsingular fibers
is an \'etale covering of $U$.

The issue is how to compactify over $B$ such a covering, in a
modular way.
We define and explore two different approaches, and compare
our resulting moduli spaces between themselves and with pre-existing 
ones.

One of the first cases in which such a question was considered is the 
one
of 2-torsion points in the Jacobian, that is, the case $N=\ma{O}$, 
$r=2$.
This was solved by A.\ Beauville \cite{beauville} by means of 
admissible
coverings.

Another remarkable, well-known instance is the one of
the so-called ``higher spin curves", which corresponds to $\N = 
\omega _f$
and, more generally, to $\N = \omega _f^{\otimes l}$.

The case of genuine ``spin curves", that is, $r=2$ (and $l=1$), was
solved by Cornalba in
\cite{cornalba1}, where a geometrically meaningful compactification is
constructed over the moduli space
$\mgbar$ of
Deligne-Mumford stable curves.
The boundary points are certain line bundles on nodal curves.

Later T.\ J.\ Jarvis (in \cite{jarvis1} and \cite{jarvis2}) used rank 
1 torsion-free sheaves to
approach the problem (over stable curves) for all $r$ and $l$.
In particular, he constructed two stack compactifications 
$\ov{\mathfrak{S}}$ and ${\text{\scshape
Root}}$, the first of which is the normalization of the second, and 
turns out to be smooth. For
$r=2$ his compactifications are all isomorphic to the one of
\cite{cornalba1}.

Other methods have recently been used by
D.\ Abramovich and T.\ J.\ Jarvis
\cite{abrjarvis}
and by A.\ Chiodo \cite{chiodo}. Both employ (differently)
``twisted curves'', which, roughly speaking, are nodal curves with a
stack structure at some nodes. The resulting compactifications are
isomorphic to Jarvis's $\ov{\mathfrak{S}}$.

Over the last decade, interest in higher spin curves has been revived 
by the
generalized Witten conjecture,
which predicts that the intersection theory on their moduli spaces is 
governed by the
Gelfand-Diki\u\i\ (also known as higher KdV) hierarchy.
This conjecture is open\footnote{Added in proof. The generalized
  Witten conjecture has recently been proved by C.\ Faber, 
S.\ Shadrin and D.\  Zvonkine in
Faber, Shadrin, Zvonkine: 
 {\it Tautological relations and the
$r$-spin Witten conjecture.}  Preprint math.AG/0612510.};
 see \cite{lee} for more details and some 
recent progress.

\bigskip
The first part of our paper presents a solution of the general 
problem in the same spirit as
\cite{cornalba1};
we complete the $r$-th root functor (and its moduli scheme)
for any family and any line bundle
$\N$ as above, by means of line bundles
on nodal curves.
For clarity (and to distinguish it from other approaches)
we name the resulting space the moduli space for {\it limit
roots}. Some of the results in this part can already be found, in one 
form
or another, in the existing literature; however, it seemed worthwhile 
to present a coherent account of the subject based on the ideas of 
\cite{cornalba1}. In order to limit the length of the
paper, we work in the context of schemes and coarse moduli spaces, 
although it is quite
clear that our methods lend themselves to the construction of stacks 
which are fine moduli
spaces for the problem at hand. We plan to return to the question in 
a future paper.

As in the case $r=2$ (see \cite{capocasa}), the combinatorial aspects 
of our construction play a
significant role, and enable us to offer a fairly explicit 
description of our compactification.

We prove in Theorem \ref{Jarvis} that, in the special case of higher 
spin curves, our moduli
space is isomorphic to the coarse space underlying Jarvis's 
$\text{\scshape Root}$.

The last section is devoted to our second approach, whose goal is to 
obtain a completion
within the compactified Picard scheme, where,
obviously, $r$-th roots on smooth curves naturally live.
Because of the relative newness of this method, we only apply it to 
higher spin curves over
$\mgbar$.
A modular compactification of the universal Picard scheme was 
constructed in \cite{caporaso} by
Caporaso, by means of line bundles on semistable curves.
This space was later given a description in terms of rank 1
torsion-free sheaves by R.\ Pandharipande \cite{pandha}, as a special 
case of a more general
result, valid for all ranks.

The question on whether the existing compactifications of higher spin 
curves
embed in the universal Picard scheme seemed, at least to us, 
irresistible.
All the more so because of the strong similarities between the 
boundary points used
by the various authors: (analogous) line bundles in \cite{cornalba1} 
and
\cite{caporaso} on the one hand, and rank 1 torsion-free sheaves in
\cite{jarvis1} and \cite{pandha} on the other.

Therefore we construct and study a new modular compactification
of higher spin curves, $\sgtil$, sitting inside the universal Picard 
scheme, and
whose main new feature is that it is not finite over $\mgbar$, as 
soon as $r\geq 3$.

We then prove that, for $r\geq 3$, the previously constructed moduli 
spaces
of limit roots, which we denote by $\sgbar$, as well as those 
constructed in \cite{jarvis1},
\cite{jarvis2}, \cite{abrjarvis}
and \cite{chiodo}, do not embed in the universal
Picard scheme.
More precisely, neither is
the natural birational map $\chi\colon \sgbar\dasharrow\sgtil$ 
regular, nor is its
inverse (Theorem \ref{comp} and its corollary). The case
$r=2$, treated by C.\ Fontanari in
\cite{fontanari}, is an exception, as $\chi$ is regular and bijective.

For higher rank, in the same spirit as ours, a compactification of 
the moduli space
of vector bundles over curves using, instead of torsion-free sheaves, 
vector bundles over
semistable curves, has recently been constructed by A.\ Schmitt
\cite{schmitt}, continuing upon
the work of C.~S.~Seshadri and others (see \cite{seshadri} and
references therein).
\subsection{Notation}
\label{not}
Throughout the paper we work over the field of complex numbers. We 
shall use the words ``line bundle'' as a synonym of ``invertible 
sheaf''.

For any scheme $Z$, we denote by $\nu=\nu_Z\colon Z^\nu\to Z$ the 
normalization morphism.

$C$ will always be a connected, nodal, projective curve of
(arithmetic) genus $g$.
We denote by
$g^{\nu}$ the
genus of its normalization.

To $C$ we associate a
{\em dual graph} $\g_C$, as follows. The
set of \emph{vertices} of $\g_C$ is the set of components of 
$C^{\nu}$. The
set of \emph{half-edges} is the set of all
points of $C^{\nu}$ mapping to a node of $C$.
Two half-edges form an \emph{edge} of $\g_C$ when the two
corresponding points map to the same node of $C$.
Recall that $g=g^{\nu}+b_1(\g_C)$, where $b_1(\g_C)$ is
the first Betti number.
Let $\Delta$ be a set of nodes of $C$, i.e., of edges of $\g_C$. We 
denote by $\ov\Delta$ the subgraph of
$\g_C$ consisting of $\Delta$ plus the abutting vertices.

We say that a nodal curve $X$ is obtained from $C$ by
{\em blowing-up} $\Delta$ if there exists $\pi\colon X\to C$ such that
$\pi^{-1}(n_i)=E_i\simeq\pr{1}$ for
any $n_i\in\Delta$, and $\pi\colon X\smallsetminus\cup_iE_i\to 
C\smallsetminus\Delta$
is an isomorphism. We call $\pi\colon X\to C$ a {\em blow-up} of $C$ 
and $E_i$ an
{\em exceptional component}.
Set $\w{X}=\ov{X\smallsetminus\cup_iE_i}$. Then $\w{\pi}\colon 
=\pi_{|\w{X}}\colon\w{X}\to C$
is the normalization of $C$ at $\Delta$.
For any $n_i\in\Delta$, set
$\{p_i,q_i\}=E_i\cap\w{X}=\w{\pi}^{-1}(n_i)$; the points $p_i$ and 
$q_i$ are called
{\it exceptional nodes}.

For families of curves,
we use the same letters but different styles
to denote the total space and a special fiber, and a star superscript 
to indicate restriction to
the complement of the special fiber. For instance, $\ma{C}\to B$ will 
denote a
family of curves with special fiber $C$, while $\ma{C}^*$ will stand
for $\ma{C}\smallsetminus C$.

Let $\ma{C}\to B$ be a family of nodal curves. We denote by $C_b$ the
fiber over $b\in B$. We say that a family of nodal curves $\X\to B$,
endowed with a $B$-morphism $\pi\colon \X\to\ma{C}$, is a
\emph{family of blow-ups} of $\ma{C}$, if for any $b\in B$ the 
restriction
$\pi_{|X_b}\colon X_b\to C_b$ is a blow-up of $C_b$.
If $\ma{N}\in\Pic{\ma C}$, we set $N_b:=\ma{N}_{|C_b}$ and, as above,
$\ma{N}^*=\ma{N}_{|\ma{C}^*}$.

We denote by $S$ a smooth, connected 1-dimensional scheme, not 
necessarily
complete, with a fixed point $s_0\in S$, and we set
$S^*=S\smallsetminus\{s_0\}$.

We shall attach to $X$ a second graph $\Sigma _X$,
whose vertices are the connected components
of $\w{X}$ and whose edges are the exceptional components of
$X$. Notice that $\Sigma _X$ is
obtained from $\g_C$ by contracting all edges
corresponding to nodes which are not blown up in $X$.

For any graph $\g$ and commutative group $G$, we denote by
$\ma{C}^0(\g,G)$ and $\ma{C}^1(\g,G)$
the groups of formal linear combinations
respectively of vertices and edges of $\g$ with
coefficients in $G$.
When we fix an orientation for $\g$,
the boundary and coboundary operators
$\partial\colon\ma{C}^1(\g,G)\to \ma{C}^0(\g,G)$
and $\delta\colon\ma{C}^0(\g,G)\to \ma{C}^1(\g,G)$ are defined in the
usual way. Then $H_1(\g,G):=\ker\partial\subset\ma{C}^1(\g,G)$.

For any positive integer $r$, $\mu_r$ is group of $r$-th roots of 
unity.
\section{Limit roots of line bundles}
Let $C$ be a nodal curve, $r$ a positive integer, and $N$ a line
bundle on $C$ whose degree is divisible by
$r$.
\subsection{The main definition}
\label{intro}
\begin{defi}
Consider a triple $(X,L,\alpha )$, where $X$
is a blow-up of $C$, $L$ is a line
bundle on $X$, and $\alpha$ is a homomorphism $\alpha\colon L^{\otimes
r}\to \pi^*(N)$ ($\pi\colon X\to C$ the natural map). We say
that $(X,L,\alpha )$ is a \emph{limit $r$-th root} of $(C,N)$ if the
following properties are satisfied:
\begin{enumerate}[(i)]
\item the restriction of $L$ to every exceptional component of $X$ has
degree 1;
\item the map $\alpha$ is an isomorphism at all points of
$X$ not belonging to an exceptional component;
\item for every
exceptional component $E_i$ of $X$, the orders of vanishing of
$\alpha$ at $p_i$ and $q_i$ add up to
$r$.
\end{enumerate}
\end{defi}
The following is an analogue,
for roots of line bundles, of the stable reduction
theorem:
\begin{prop}
\label{ssred}
Let $\mathcal{C}\to S$ be a family of nodal curves,
smooth over $S^*$.
Consider line bundles
$\ma{N}\in\Pic\ma{C}$ and $\ma{L}\in\Pic\ma{C}^*$ with an
isomorphism $\alpha\colon\ma{L}^{\otimes
r}\stackrel{\sim}{\to}
\mathcal{N}_{|\ma{C}^*}$.
Then there exist:
\begin{enumerate}[1)]
\item A diagram
\hfill
{}{}$\xymatrix{{\mathcal{X}}\ar@/^1pc/[rr]^{h}\ar[r]\ar[dr]&
{}{}{\mathcal{C}'}\ar[r]\ar[d]&{\mathcal{C}}\ar[d]\\&{S'}\ar[r]&{S}
{}{} }$\hskip40pt\hfill\ \\

where $S'$ is a smooth curve; $S'\to S$ is a finite morphism of degree
$r$, \'etale over $S^*$; $\mathcal{C}':=\mathcal{C}\times _S S'$;
$\mathcal{X}\to\mathcal{C}'$ is an
isomorphism outside the central fiber; $X$ is a blow-up of $C'$.
\item A line bundle $\ma{L}'\in\Pic\mathcal{X}$ such that
$\ma{L}'_{|\ma{X}^*}=h^*(\ma{L})$.
\item An isomorphism
$\alpha'\colon(\ma{L}')^{\otimes r}\stackrel{\sim}{\to}
h^*(\ma{N})(-D)$, where $D$
is an effective Cartier divisor supported on the exceptional
components of the central fiber,
such that $\alpha'_{|\ma{X}^*}=h^*(\alpha)$ and
$(X,\ma{L}'_{|X},\alpha'_{|X})$ is a limit $r$-th root of
$(C,\ma{N}_{|C})$.
\end{enumerate}
Moreover, for any $(\X_1\to\ma{C}',\ma{L}'_1,\alpha'_1,D_1)$ and
$(\X_2\to\ma{C}',\ma{L}'_2,\alpha'_2,D_2)$ satisfying the above
conditions, there exists an isomorphism
$\sigma\colon\X_1\stackrel{\sim}{\to}\X_2$ over
$\ma{C}'$ such that $\sigma^*(D_2)=D_1$, $\sigma^*(\ma{L}_2')\simeq
\ma{L}_1'$, and this last isomorphism is compatible with $\alpha'_1$
and $\alpha'_2$.
\end{prop}
The proof is in \cite{jarvis1}, section \S 4.2.2.
For more details on the work of Jarvis
on the subject, see section~\ref{spin}.

The notion of limit root generalizes to families in a natural
way. Let $\ma{C}\to T$ be a family of nodal curves and
$\ma{N}\in\Pic\ma{C}$. A \emph{limit $r$-th root} $(\ma{X}, \ma{L},
\alpha)$ of $(\ma{C},\ma{N})$ is the datum of
a family $\pi\colon \ma{X}\to \ma{C}$
of blow-ups of $\ma{C}$,
a line bundle $\ma{L}\in\Pic \ma{X}$ and a homomorphism
$\alpha\colon\ma{L}^{\otimes r}\to\pi^*(\ma{N})$, such that for
all $t\in T$, $(X_t,L_t,\alpha_t)$ is a limit $r$-th root of
$(C_t,N_t)$.

We can restate
Proposition \ref{ssred} by saying that if $\ma{C}\to S$ is a
family of nodal
curves, $\ma{N}\in\Pic\ma{C}$ and $(\ma{C}^*,\ma{L},\alpha)$
is an $r$-th root of
$\ma{N}_{|\ma{C}^*}$, then for some finite covering
$S'\to S$ of degree $r$ there exists a limit $r$-th root
$(\ma{X},\ma{L}',\alpha')$ of the pull-back of $\ma{N}$ to the
base-changed family $\ma{C}' \to S'$ which extends over all of $S'$ 
the
pull-back of $(\ma{C}^*\to S^*,\ma{L},\alpha)$.
\subsection{Weighted graphs}
\label{wg}
To a limit $r$-th root
$(X,L,\alpha)$ we shall attach a {\it weighted subgraph} of the dual
graph $\g_C$ of $C$. The underlying graph is the subgraph $\ov\Delta$ 
corresponding to the set
$\Delta$ of nodes which are blown up in $\pi\colon X\to C$. The 
weight function $w$ is defined on
the exceptional nodes (or, equivalently, on the half-edges of 
$\ov\Delta$),
takes values in
$\{1,\dots,r-1\}$, and is given by
$w(p_i)=u_i$,
$w(q_i)=v_i$, where
$u_i,v_i$ are the orders of vanishing of $\alpha$ at the points 
$p_i,q_i$ respectively. We denote
this weighted graph by $\Delta^w$.
Clearly, the following
conditions are satisfied:
\begin{enumerate}[\bf (C1)]
\item $u_i+v_i=r$ for any $i$;
\item
for every irreducible component $C_j$ of $C$,
the sum of all weights assigned to
the vertex corresponding $C_j$ is congruent to $\deg_{C_j} N$ modulo
$r$.
\end{enumerate}
Conversely, for any assignment of a
weighted subgraph $\Delta^w$ of
$\g_C$ satisfying conditions (C1) and (C2), there
is a limit $r$-th root of $N$ whose graph is
$\Delta^w$. In fact, let
$\w{\pi}\colon \w{X}\to C$ be the partial
normalization of $C$ at $\Delta$.
By (C2),
the line bundle $\w{\pi}^*(N)(-\sum_{n_i\in\ev}
(u_ip_i+v_iq_i))$ on $\w{X}$ admits an $r$-th root $\w{L}
\in\Pic\w{X}$. Let $X$ be the blow-up of $C$ at $\Delta$.
We obtain a limit root $(X,L,\alpha )$ of $N$
by gluing to $\w{L}$ (however we wish)
the degree one line bundle on each exceptional component; the map 
$\alpha\colon
L^{\otimes r}\to\pi^*(N)$ is defined
to agree
with $(\w{L})^{\otimes 
r}=\w{\pi}^*(N)(-\sum_i(u_ip_i+v_iq_i))\hookrightarrow
\w{\pi}^*(N)$ on $\w{X}$, and to be zero on all exceptional
components.
\begin{remark}
\label{graphs}
The {\it multidegree} $\un{d}=\un{\deg}\ N$ of $N$ is the set of the 
degrees
$\deg_{C_j} N$ for all irreducible components $C_j$ of $C$. Its class
modulo $r$ can be viewed as an element
$\dnmodr$ of $ \ma{C}^0(\g_C,\Z/r)$ defined naturally as
$$
\dnmodr:=\sum_j \Bigl(\deg_{C_j}N \pmod r\Bigr) [C_j].$$
Of course, condition (C2) above only depends on $\dnmodr$.

Now fix an orientation for $\g_C$;
then there is a natural bijection between $\ma{C}^1(\g_C,\Z/r)$
and the set of
weighted subgraphs of $\g_C$ satisfying (C1): to a weighted subgraph 
$\ev^w$ we
associate the 1-chain $\sum_i u_i[n_i]$, with the convention that
the edge $n_i$ is oriented entering the vertex corresponding to
$u_i$.

In this set-up, the boundary operator
$\partial\colon\ma{C}^1(\g_C,\Z/r)\to \ma{C}^0(\g_C,\Z/r)$
can be viewed as a map from weighted subgraphs of $\g_C$ satisfying 
(C1) to
multidegrees reduced modulo $r$.
Then
$\ev^w$ satisfies (C2)
with respect to $\dnmodr$
if and only if
$\partial(\ev^w)=\dnmodr$.
Since the cardinality of $\partial^{-1}(\dnmodr)$ equals
that of $\ker\partial=H_1(\g_C,\Z/r)$, we see that:
\begin{enumerate}[$\bullet$]
\item
{\em for every $N$ as above, there exist
$r^{b_1(\g_C)}$ weighted subgraphs of $\g_C$ satisfying (C1) and (C2)
with respect to $N$}.
\end{enumerate}
\end{remark}
\subsection{Automorphisms of limit roots}
\label{auto}
\begin{defi}
\label{iso}
An \emph{isomorphism} of limit $r$-th roots $(\ma{X}\to
T,\ma{L},\alpha)$ and $(\ma{X}'\to T,\ma{L}',\alpha')$ of 
$(\ma{C},\ma{N})$
is the datum of:
\begin{enumerate}[(a)]
\item an isomorphism
$\sigma\colon\ma{X}\to\ma{X}'$ over $\ma{C}$
\item an isomorphism
$\tau\colon\sigma^*\ma{L}'\to\ma{L}$ that makes the following diagram
commute:
\[
{}{}\xymatrix{{\sigma^*(\mathcal{L}')^{\otimes
{}{} r}}\ar[d]_{\sigma^*(\alpha')}\ar[r]^{\phantom{a}\tau^{\otimes r}}
{}{}&{\mathcal{L}^{\otimes r}}\ar[d]^{\alpha}\\
{}{}{\sigma^*(\pi')^*(\ma{N})}\ar[r]^{\phantom{aa}\sim}&{\pi^*(\ma{N})} 
}\hskip25pt
\]
\end{enumerate}
\end{defi}
This notion of isomorphism agrees with the one of \cite{cornalba2}, 
rather than with the
one of \cite{cornalba1}. At the level of automorphisms, for instance, 
the difference is that
here fiber multiplication in $\mathcal{L}$ by a non-trivial root of 
unity is viewed as a
non-trivial automorphism of the limit root.
To make contact with the terminology of \cite{cornalba1} and 
\cite{cornalba2}, an
isomorphism between $\X$ and $\X'$ over $\ma{C}$ will sometimes be 
referred to as an
\emph{inessential} isomorphism.
Likewise, automorphisms of $\X$ over $\ma{C}$ will sometimes be 
called inessential
automorphisms of $\X$. As customary, the group of these automorphism 
will be denoted
$\Aut_{\ma{C}}(\X)$. Given a limit root $(\X\to T,\ma{L}, \alpha)$,
we denote by $\Aut(\X\to T,\ma{L}, \alpha)$ the group of its 
automorphisms.
Clearly, $\Aut(\X\to T,\ma{L}, \alpha)$ maps to $\Aut_{\ma{C}}(\X)$; 
the kernel is naturally isomorphic to a product of copies of $\mu_r$, 
one for each connected component of $T$.

We now state and prove some of the basic properties of automorphisms 
of limit roots. Let $X$ be a blow-up of $C$, and let $E_1,\dots,E_m$ 
be its exceptional components. View each $E_i$ as a copy of the 
Riemann sphere, with the origin placed at $q_i$ and the point at 
infinity placed at $p_i$. Then any inessential automorphism of $X$ 
acts on each $E_i$ as multiplication by a nonzero constant $l_i$, and 
any assignment of the $l_i$ corresponds to an inessential 
automorphism. This proves the first part of the following result.
\begin{lemma}
\label{iness}
Let $(X,L,\alpha)$ be a limit root of $(C,N)$, and fix an orientation 
on the
graph $\Sigma_X$. Then:
\begin{enumerate}[(i)]
\item There is a natural identification $\Aut_C(X)\simeq 
\ma{C}^1(\Sigma _X,\C^*)$.
\item There is a natural identification 
$\Aut(X,L,\alpha)\simeq\ma{C}^0(\Sigma _X,\mu_r)$, and the 
homomorphism $\Aut(X,L,\alpha)\to\Aut_C(X)$ corresponds
to the composition of the coboundary map $\delta\colon\ma{C}^0(\Sigma 
_X,\mu_r)\to
\ma{C}^1(\Sigma _X,\mu_r)$ with the inclusion $\ma{C}^1(\Sigma
_X,\mu_r)\to\ma{C}^1(\Sigma _X,\C^*)$.
\end{enumerate}
\end{lemma}
\begin{proof}We have already proved {\it (i)}. To prove {\it (ii)}, 
it is convenient to work with the geometric bundles $\mathbb{V}(L)$ 
and $\mathbb{V}(L^{\otimes r})$ rather than with $L$ and $L^{\otimes 
r}$. In terms of these, an automorphism of $(X,L,\alpha)$ is a pair 
$(\sigma,\tau)$, where $\sigma$ is an inessential automorphism of $X$ 
and $\tau$ is an automorphism of $\mathbb{V}(L)$ which is compatible 
with $\sigma$, is linear on the fibers of $\mathbb{V}(L)\to X$, and 
is such that $\alpha\circ\tau^{\otimes r}=\alpha$. This last 
condition means that, above each connected component $X_j$ of 
$\widetilde{X}$, $\tau^{\otimes r}$ is the identity, and hence $\tau$ 
must be the fiberwise multiplication by an $r$-th root of unity 
$l_j$. If $E_i$ is an exceptional component connecting $X_j$ and 
$X_k$, the restrictions of $\tau$ and $\sigma$ to 
$\mathbb{V}(L)_{|E_i}\to E_i$ can thus be viewed as an automorphism 
pair of $\mathbb{V}(\ma{O}(1))\to\mathbb{P}^1$ which is the 
multiplication by $l_j$ on the fiber above the origin and the 
multiplication by $l_k$ on the fiber above the point at infinity. 
These pairs are easy to describe. View $\mathbb{P}^1$ as the 
projective space $\mathbb{P}V$ constructed on a two-dimensional 
vector space $V$, and its origin and point at infinity as $[v_0]$ and 
$[v_1]$, where $v_0,v_1$ is a basis of $V$; then 
$\mathbb{V}(\ma{O}(1))$ is the set of all pairs $(\ell,\ph)$, where 
$\ell$ is a one-dimensional subspace of $V$ and $\ph$ an element of 
its dual. An automorphism pair as above then acts on 
$\mathbb{V}(\ma{O}(1))$ via $(\ell,\ph)\mapsto (f(\ell),\ph\circ 
f^{-1})$, where $f$ is the automorphism of $V$ such that 
$f(v_{0})=l_j^{-1}v_{0}$, $f(v_{1})=l_k^{-1}v_{0}$, and on 
$\mathbb{P}V$ via $\ell\mapsto f(\ell)$. Thus $\sigma$ acts on $E_i$ 
as multiplication by $l_j/l_k$. This means that giving an 
automorphism of $(X,L,\alpha)$ is the same as giving an $r$-th root 
of unity for each connected component of $\widetilde{X}$, that is, an 
element of $\ma{C}^0(\Sigma _X,\mu_r)$, and that the map 
$\Aut(X,L,\alpha)\to\Aut_C(X)$ is as claimed.
\end{proof}
It follows from \ref{iness} that $\Aut(X,L,\alpha)$
has cardinality $r^{\gamma}$, where $\gamma$ is
the number of connected components of $\w{X}$.

\medskip
Fix $(C,N)$ 
and $\pi:X\to C$ as usual.
We need to further investigate the action 
of $\Aut_C(X)$ on certain line
bundles on $X$.
Consider
pairs 
$(M,\beta)$,
where $M$ is a line bundle on $X$,
$\beta\colon 
M\to\pi^*(N)$ is a homomorphism which is an isomorphism
on the 
complement of the exceptional components of $X$,
and $M$ is assumed 
to have degree $r$ on each
exceptional component. This degree 
condition forces
$\beta$ to vanish identically on each exceptional 
component

An example of this situation is provided by
$\alpha\colon 
L^{\otimes r}\to
\pi^*(N)$, where $(X,L,\alpha)$ is a limit $r$-th 
root.
\begin{lemma}
\label{pull}
Let $(M,\beta)$ and $(M',\beta')$ be two pairs as above. 
Then:
\begin{enumerate}[(i)]
\item There exist $\sigma\in\Aut_C(X)$ and an isomorphism 
$\vartheta\colon M\to \sigma^*(M')$ such that 
$\beta=\sigma^*(\beta')\circ\vartheta$ if and only if, for any 
exceptional node $q$, the orders of vanishing of $\beta$ and $\beta'$ 
at $q$ are the same.
\item Under the identification (i) of 
\ref{iness}, the set of those $\sigma\in \Aut_C(X)$ such that there 
is an isomorphism $\vartheta\colon M\to \sigma^*(M)$ with the 
property that $\beta=\sigma^*(\beta)\circ\vartheta$ corresponds to 
$\ma{C}^1(\Sigma _X,\mu_r)$.
\end{enumerate}
\end{lemma}
\begin{proof}
The proof is similar to the one of \ref{iness}. The 
``only if'' part of {\it (i)} is clear. To prove the converse, denote 
by $u_i$ (resp., $v_i$) the order of vanishing of $\beta$ and 
$\beta'$ at $p_i$ (resp., $q_i$). The homomorphisms $\beta$ and 
$\beta'$ determine isomorphisms $M_{|\w{X}}\simeq\pi^*(N)(-\sum 
(u_ip_i+v_iq_i)) \simeq M'_{|\w{X}}$, hence there is a unique 
isomorphism $\vartheta'\colon M_{|\w{X}}\to M'_{|\w{X}}$ such that 
$\beta_{|\w{X}}=\beta'_{|\w{X}}\circ\vartheta'$. Clearly, $\vartheta$ 
must restrict to $\vartheta'$ on $\w{X}$. It remains to construct 
$\sigma$ and $\vartheta$ on the exceptional components. If $E_i$ is 
one of these, we identify the restrictions to $E_i$ of 
$\mathbb{V}(M)\to X$ and $\mathbb{V}(M')\to X$ to 
$\mathbb{V}(\ma{O}(r))\to\mathbb{P}(V)$, where $V$ is a 
two-dimensional vector space. We also think of $\mathbb{P}(V)$ as 
being attached to the rest of $X$ at $[v_0]$ and $[v_1]$, where 
$v_0,v_1$ is a basis of $V$. The isomorphisms $\vartheta$ and 
$\sigma$, if they exist, restrict on 
$\mathbb{V}(\ma{O}(r))\to\mathbb{P}(V)$ to an automorphism pair 
fixing $[v_0]$ and $[v_1]$. Now, $\mathbb{V}(\ma{O}(r))$ is the set 
of all pairs $(\ell,\ph^{\otimes r})$, where $\ell$ is a 
one-dimensional subspace of $V$ and $\ph$ an element of its dual. Any 
automorphism pair fixing $[v_0]$ and $[v_1]$ then acts on 
$\mathbb{V}(\ma{O}(r))$ via $(\ell,\ph^{\otimes r})\mapsto 
(f(\ell),(\ph\circ f^{-1})^{\otimes r})$, where $f$ is an 
automorphism of $V$ such that $f(v_{0})=b_0v_{0}$, 
$f(v_{1})=b_1v_{0}$, $b_0,b_1\in\mathbb{C}^*$, and on $\mathbb{P}V$ 
via $\ell\mapsto f(\ell)$. The action of such an automorphism pair on 
the fiber above $[v_i]$, $i=0,1$, corresponds to multiplication by 
$b_i^{-r}$, hence we can choose $b_0$ and $b_1$ so that it matches 
the one of $\vartheta'$; the choices are unique up to multiplication 
by $r$-th roots of unity. This proves the existence of $\sigma$ and 
$\vartheta$. It also proves {\it (ii)}. In fact, $\sigma$ acts on 
$E_i$ as multiplication by $b_1/b_0$. In the setup of {\it (ii)}, 
both $b_0$ and $b_1$ are (arbitrary) roots of unity. The conclusion 
is that $\sigma$ acts on each exceptional component as multiplication 
by an $r$-th root of unity, and conversely that any element of 
$\ma{C}^1(\Sigma _X,\mu_r)$ comes from an automorphism $\sigma$ as in 
{\it (ii)}.
\end{proof}
\subsection{The general moduli problem}
\label{moduli}
Let $B$ be a scheme,
$f\colon\ma{C}\to B$ a family of nodal curves, $\ma{N}\in
\Pic\ma{C}$ a
line bundle of relative degree divisible by $r$. We shall use the
following notation: for any $B$-scheme $T$ with structure
map $p\colon T\to B$ we set
$$
{}{}\xymatrix{
{}{}{\ma{C}_T:= \ma{C}\times _B 
T}\ar[r]^(.7){\ov{p}}\ar[d]&{\ma{C}}\ar[d]^f
{}{}\\ {T}\ar[r]^{p}& B}
$$
where $\overline{p}\colon\ma{C}_T\to \ma{C}$ is the projection. The
pull-back of $\ma{N}$ to $\ma{C}_T$ is
denoted by
$$\ma{N}_T:=\overline{p}^*\ma{N}.$$
Having fixed
$B$, $f\colon\ma{C}\to B$ and $\ma{N}$ as above, we introduce
a contravariant functor
$$\snf \colon \{B{\text {-schemes}}\}\la
\{{\rm {sets}}\},$$
encoding the moduli problem for
limit $r$-th roots of $\ma{N}$, as follows. For every $B$-scheme
$T$, the set $\snf (T)$ is the set of all
limit $r$-th roots of $\ma{N}_T$ (which is, with the above
notation, a line bundle on the total space of the
family $\ma{C}_T\to T$), modulo isomorphisms of limit roots.

For every morphism of $B$-schemes $\rho\colon T\to T'$, the
corresponding map
$$\snf(\rho)\colon\snf (T')\la \snf (T)$$
sends a
limit $r$-th root $(\ma{X}',\ma{L}', \alpha ')$ of
$\ma{N}_{T'}$, to the limit root $(\ma{X},\ma{L}, \alpha)$ of
$\ma{N}_T$ defined by $\ma{X}:=\ma{X}'\times _{T'}T$
and $\ma{L}:=\ov{\rho}^*\ma{L}'$ where $\ov{\rho}\colon\ma{X}\to 
\ma{X}'$
is the projection. Finally, the morphism $\alpha '$ is the
natural pull-back of $\alpha$.

We shall prove
\begin{thm}
\label{coarsemoduli}
Let $f\colon\ma{C}\to B$ be a
family of nodal curves, with $B$ a quasi-projective
scheme. Let $\ma{N}\in \Pic\ma{C}$ be a line bundle on $\ma{C}$ of
relative degree divisible by $r$.

Then the functor $\snf$ is
coarsely represented by a quasi-projective scheme $\snbar$, finite 
over $B$.
If $B$ is projective, then $\snbar$
is projective.
\end{thm}
We shall denote by $\sn$ the open subscheme of $\snbar$ parametrizing
(closed)
points
whose underlying curve has no exceptional components.
\begin{prop}
\label{smooth}
Let $\xi=(X,L,\alpha)\in\snbar$. The morphism
$p\colon\snbar\to B$ is smooth at $\xi$ if and only if
$b_1(\Sigma_X)=0$.

In particular, if $B$ is smooth,
$\snbar$ is smooth at all points $\xi$ such that $X$
has no exceptional components, or is of compact type.
\end{prop}
The proofs of these results are in section \ref{proof}.
\begin{remark}
\label{remark}
A purely formal consequence of Theorem \ref{coarsemoduli}
is the following useful fact.
Let $B'$ be a $B$-scheme, and denote
by
$$f'\colon\ma{C}'=\ma{C}\times_BB'\la B'$$
the family obtained by
base change, and by $\ma{N}'$ the pull-back of
$\ma{N}$ to $\ma{C}'$. Then we naturally have
$$\ov{S}^{\,r}_{f'}(\mathcal{N'})\simeq
\snbar\times_B B'.$$
\end{remark}
\begin{remark}
Let $f\colon\ma{C}\to B$ and $\ma{N}\in\Pic\ma{C}$ be as above. For
any $\ma{L}_0\in\Pic\ma{C}$, the functors $\snf$ and
$\ov{\mathcal{S}}^{\,r}_f(\mathcal{N}\otimes\ma{L}_0^{\otimes r})$ are
isomorphic, so the same holds for the corresponding moduli spaces
$\snbar$ and $\ov{S}^{\,r}_f(\ma{N}\otimes\ma{L}_0^{\otimes r})$.
This is because for any limit $r$-th root $(\ma{X},\ma{L})$ of
$\ma{N}$, $(\ma{X},\ma{L}\otimes\pi^*(\ma{L}_0))$ is a limit $r$-th 
root of
$\ma{N}\otimes\ma{L}_0^{\otimes r}$.
\end{remark}
\begin{remark}
\label{rel}
Let $f\colon\ma{C}\to B$ and $\ma{N}\in\Pic\ma{C}$ be as above.
For any positive integer $s$, there is a
natural injective morphism
\[J\colon\snbar\la\ov{S}^{\,rs}_f(\ma{N}^{\otimes s})\]
over $B$.
Indeed for any limit $r$-th root $(\ma{X},\ma{L},\alpha)$ of
$\ma{N}$, $(\ma{X},\ma{L},\alpha^{\otimes s})$ is a limit $(rs)$-th 
root of
$\ma{N}^{\otimes s}$. The image of $J$ is a union of irreducible
components of $\ov{S}^{\,rs}_f(\ma{N}^{\otimes s})$.
\end{remark}
\section{Proof of the existence theorem}
\label{dimo}
We shall now construct $\snbar$, and prove Theorem 
\ref{coarsemoduli}. First of all, working in the analytic category, 
we construct the universal deformation of a limit root. Then we glue 
together the bases of the universal deformations, in order to give a 
complex algebraic structure to $\snbar$.

We will make use of the following simple results. The first one is 
well-known, while the second is Lemma 1.1 of \cite{cornalba1}.
\begin{remark}
\label{fibrati}
Let $\ma{Y}\to T$ be a family of nodal curves, $t_0$ a point of $T$, 
and $Y$
the fiber over $t_0$. Consider two line bundles
$\ma{P}\in\Pic\ma{Y}$ and $L\in\Pic Y$ such that there exists an 
isomorphism $\iota_0\colon L^{\otimes r}\to\ma{P}_{|Y}$. Then,
up to shrinking $T$, there exist a line bundle $\ma{L}\in\Pic\ma{Y}$
extending $L$ and an isomorphism $\iota\colon \ma{L}^{\otimes 
r}\to\ma{P}$ extending
$\iota_0$. Moreover, this extension is unique up to isomorphism, 
meaning that if $(\ma{L}',\iota')$ is another extension, then there 
is an isomorphism
$\chi\colon \ma{L}\to\ma{L}'$ on a neighbourhood of $Y$ which 
restricts to the identity on
$L$ and is such that $\iota=\iota'\circ\chi^{\otimes r}$.
\end{remark}
\begin{remark}
\label{lemmino}
Let $\ma{Y}\to T$ be a 
family of nodal curves, and let $E$ be a smooth rational component of 
one of its fibers. If $\ma{L}$ is a line bundle on $\ma{Y}$ whose 
restriction to $E$ has zero degree, then $\ma{L}$ is trivial on a 
neighbourhood of $E$ in $\ma{Y}$.
\end{remark}
\subsection{Local structure near exceptional curves}
\label{def}
Let $\ma{C}\to T$ be a family of nodal curves, $\ma{N}\in\Pic\ma{C}$
and
$(\ma{X}\to T,\ma{L},\alpha)$ a limit $r$-th
root of $(\ma{C},\ma{N})$. Fix an exceptional component
$E\subset\ma{X}$ lying over a point $t_0\in T$, let $u,v$ be the 
corresponding weights,
and denote by $n$ the point of $\ma{C}$ to which $E$ contracts. As 
usual,
$\pi\colon \ma{X}\to\ma{C}$ is the natural map; recall that
$\pi^*\ma{N}$ is trivial on a neighbourhood of $E$. We shall describe 
the geometry of
$\X \to T$ near $E$.

The prototype is the following family:
for the base, consider the affine plane curve
$C_{u,v}\subset\mathbb{A}^2_{\,w,z}$ \footnote{The notation
$\mathbb{A}^2_{\,w,z}$ means, here and similarly below, that
$w,z$ are the coordinate
functions on $\mathbb{A}^2$.} given by the equation
$w^{u}=z^{v}$.
The total space is the surface
\[S_{u,v}=\{xs_0=w s_1,\ ys_1=z s_0,\ w^{u}=z^{v}\}
\subset \mathbb{A}^4_{\,x,y,w,z}\times\pr{1}_{s_0:s_1}.\]
The projection to $\mathbb{A}^2_{\,w,z}$
yields a family of open curves
$$
f_{u,v}\colon S_{u,v} \la C_{u,v}
$$
whose only singular fiber lies over the
point $(0,0)\in C_{u,v}$ and contains an exceptional component, i.e., 
the locus of all points of $S_{u,v}$ for which $x=y=w=z=0$.

Notice that $S_{u,v}$ has $d=(u,v)$ irreducible components, each one
isomorphic to $S_{u/d,\,v/d}$.
For instance, when $d=1$, $S_{u,v}$ is smooth if and only if $u=v=1$ 
(in
which case it is a resolution of an $A_1$-singularity); it is normal
if and only if $u=1$ or $v=1$ (in which case it has an isolated
$A_{v-1}$ or $A_{u-1}$ singularity in one exceptional node). If $u>1$
and $v>1$, $S_{u,v}$ is singular along the whole special fiber of 
$f_{u,v}$
(the curve $w=z=0$).

We now show that the morphism $f_{u,v}\colon S_{u,v} \to C_{u,v}\ $ 
is a model for $\X\to
T$ near $E$. To set up, recall that a neighbourhood $V$ of $n$ in 
$\ma{C}$ is of the form
$\{xy=\ell\}$, where $\ell$ is a regular function on a neighbourhood 
$T_0$ of $t_0$. In
other terms, there is a fiber product diagram
\begin{equation}
\label{liftandus}
{}{}\xymatrix{
{}{}\ma{C} \supset\hskip-38pt &V\ar[r]\ar[d]&V'\ar[d]^e\\
{}{}&{T_0}\ar[r]^{\ell}& \mathbb{A}^1_s}
\end{equation}
where $V'$ is a neighbourhood of the origin in $\mathbb{A}^2_{x,y}$ 
and $e$ is given by
$s=xy$.

The curve $C_{u,v}$ maps to $\mathbb{A}^1_s$ via $(w,z)\mapsto wz$, 
and $S_{u,v}$ to
$\mathbb{A}^2_{x,y}$ via projection. We will show that $\ell$ lifts 
to a morphism
$h\colon T_0\to C_{u,v}$. More precisely, we shall prove the 
following result.
\begin{lemma}
\label{local}
Possibly after shrinking $T_0$,
 diagram 
\eqref{liftandus} lifts to
a fiber product diagram
\[{}{}\xymatrix{
{}{}\X \supset\hskip-38pt &U\ar[r]\ar[d]&\text{\phantom{U}}
{}{}\ar[d]^{f_{u,v}}&\hskip-43pt U'\subset S_{u,v}\\
{}{}\ 
&{T_0}\ar[r]^{h}& C_{u,v}}
\]
where $U$ is a neighbourhood of $E$ in 
$\ma{X}$ and $U'$ is a
neighbourhood of the exceptional curve in 
$S_{u,v}$.
Moreover, the lifting can be chosen so that the following 
holds.
Let $\zeta_{\ma{N}}$ be any generator of $\pi^*\ma{N}$ over 
$U$, and
denote by
$U_0$ (resp., $U_1$) the complement in
$U$ of the 
inverse image of the locus $s_1=0$ (resp.,
$s_0=0$). Then there are 
generators
$\zeta_0$ and $\zeta_1$ of $\ma{L}^{\otimes r}$ over $U_0$ 
and $U_1$ such that
$$
s_0^r\zeta_0=s_1^r\zeta_1\text{ on } U_0\cap 
U_1,\quad
\alpha(\zeta_0)=x^u\zeta_{\ma{N}}
\text{ on }U_0,\quad 
\alpha(\zeta_1)=y^v\zeta_{\ma{N}}\text{ on }U_1.\ \footnote{For 
simplicity, we continue to use $x$ and $y$ to denote the pull-backs 
of these functions to $U$.}
$$
\end{lemma}
\begin{proof}
A weaker lifting property than the one claimed obviously holds. 
Consider the threefold
\[Z=\{xs_0=w s_1,\ ys_1=z s_0\}
\subset \mathbb{A}^4_{\,x,y,w,z}\times\pr{1}_{s_0:s_1}.\]
The projection to $\mathbb{A}^2_{\,w,z}$
yields a family of open curves
$\phi\colon Z \la \mathbb{A}^2_{\,w,z}$, whose restriction to 
$C_{u,v}$ is clearly just $S_{u,v}
\to C_{u,v}$. Then, after suitable shrinkings, diagram 
\eqref{liftandus} lifts (non uniquely)
to
\[{}{}\xymatrix{
{}{}U\ar[r]\ar[d]&W\ar[d]^{\phi}\\
{}{}{T_0}\ar[r]^{k}& \mathbb{A}^2_{\,w,z}}
\]
where $W\subset Z$ is a neighbourhood of the exceptional curve in the 
central
fiber of $\phi$. In concrete terms, this means that
$U$ has equations:
\[ x\eta_0=k_0\ \text{in
}U_0,\quad y\eta_1=k_1\ \text{in }U_1, \]
where $\eta_0\eta_1=1$ on $U_0\cap U_1$, $x$ and $y$ vanish on $E$,
$\eta_0$ and $\eta_1$ vanish on $\w{X}\cap U$,
$k_0$ and $k_1$ are regular functions on $T_0$ which vanish at $t_0$. 
Clearly, $\eta_0$
is the pull-back of the quotient $s_0/s_1$, and $\eta_1$ is its 
inverse. What must be
shown is that $k$ can be chosen in such a way that $k_0^u=k_1^v$. 
This is where the
presence of a limit root comes into play.

Consider the line bundle $\ma{M}$ on $U$ whose transition function, 
relative
to the cover $\{U_0,U_1\}$, is $\eta_1^r\in\ma{O}^*(U_0\cap U_1)$. 
Clearly, $\deg
\ma{M}_{|E}=\deg\ma{L}^{\otimes r}_{|E}=r$. By \ref{lemmino}, we may 
assume that $\ma{L}^{\otimes r}_{|U}\simeq\ma{M}$, possibly after 
shrinking $U$. Thus there are
generators $\zeta'_0$, $\zeta'_1$ of $\ma{L}^{\otimes r}$, over $U_0$ 
and $U_1$
respectively, such that $\zeta'_0=\eta_1^r\zeta'_1$ on $U_0\cap U_1$.

Let $\zeta_{\ma{N}}$ be a generator of $\pi^*(\ma{N})$ over $U$. We 
have
$\alpha(\zeta'_i)=c_i \zeta_{\ma{N}}$ on $U_i$, $i=0,1$, where $c_0$ 
and $c_1$ are
regular functions on $U_0$ and $U_1$ respectively.
Consider the expressions of $c_0$ and $c_1$ as power
series:
\begin{equation}
\label{develop}
c_0=\sum_{i\geq 0}a_ix^i+\sum_{i< 0}a_i\eta_0^{-i}
\ \text{ in $U_0$,}\qquad
c_1=\sum_{i\geq 0}b_iy^i+\sum_{i< 0}b_i\eta_1^{-i}
\ \text{ in $U_1$},
\end{equation}
where the $a_i$ and the $b_i$ are regular functions on $T_0$.
Since $\pi^*(\ma{N})$ has degree zero on $E$ while $\ma{L}^{\otimes 
r}$ has
degree $r$, $\alpha$ is identically zero on $E$, which gives
$a_i(t_0)=b_i(t_0)=0$ for all $i\leq 0$. On the other hand, by 
definition, the order of
vanishing of $(c_0)_{|X\cap U_0}$ at $p$ is $u$ and that of
$(c_1)_{|X\cap U_1}$ at $q$ is $v$, so $a_u(t_0)\neq 0$ and 
$b_v(t_0)\neq 0$. Up to
shrinking $T_0$, we may thus assume that $a_u$ and $b_v$ are units.

On $U_0\cap U_1$ we have
$\alpha(\zeta'_0)=\eta_1^r\alpha(\zeta'_1)$, so $c_0=\eta_1^r
c_1$. We also have $\eta_0=\eta_1^{-1}$, $x=k_0\eta_1$, 
$y=k_1\eta_1^{-1}$.
Substituting, we get:
\begin{equation}
\label{compare}
\sum_{i\geq 0}a_ik_0^i\eta_1^i+\sum_{i< 0}a_i\eta_1^{i}
=\sum_{i\geq 0}b_ik_1^i\eta_1^{r-i}+\sum_{i< 0}b_i\eta_1^{r-i}
\ \text{ in $U_0\cap U_1$}.
\end{equation}
Comparing terms of degree $u$ in $\eta_1$, we get
$a_uk_0^u=b_vk_1^v$. This is not exactly
what we are asking for, but we may proceed as follows. Let $\delta$ 
be an $r$-th root of
$b_v/a_u$ (a unit), and define a morphism $h\colon T_0\to 
\mathbb{A}^2_{w,z}$ by setting
$h_0=\delta^{-1} k_0$,
$h_1=\delta k_1$. Clearly, $h_0^u=h_1^v$, and $h$ is a lifting of 
$\ell$, since
$h_0h_1=k_0k_1=\ell$. A compatible lifting of $V\to 
\mathbb{A}^2_{x,y}$ to a morphism $U\to S_{u,v}$
is given by
\begin{center}
$(x,\eta_0,t)\mapsto 
(x,\eta_0k_1(t),h_0(t),h_1(t),[\delta^{-1}\eta_0:1])
\quad\text{on }U_0$,\par
$(y,\eta_1,t)\mapsto (\eta_1k_0(t),y,h_0(t),h_1(t),[1:\delta\eta_1])
\quad\text{on }U_1$,
\end{center}
where $t$ stands for a variable point in $T_0$.

It remains to find the local generators $\zeta_0$ and $\zeta_1$.
What we have done so far amounts to showing that, replacing
$\eta_0$ with $\delta^{-1}\eta_0$ and $\eta_1$ with $\delta\eta_1$, 
we may assume
that $U$ is defined by equations $x\eta_0=h_0$ on $U_0$ and 
$y\eta_1=h_1$ on $U_1$,
where $\eta_0\eta_1=1$ on $U_0\cap U_1$ and $h_0^u=h_1^v$. Moreover, 
replacing
$\zeta_0'$ with $\delta^u\zeta_0'$ and $\zeta_1'$ with 
$\delta^{-v}\zeta_1'$, we may
also assume that the coefficients $a_u$ and $b_v$ in the power series 
developments
\eqref{develop} are equal, while conserving the property that
$\zeta_0'=\eta_1^r\zeta_1'$. Now, comparing terms of all degrees in 
identity
\eqref{compare}, we conclude that
$$
c_0=x^u\rho,\quad c_1=y^v\rho,
$$
where
$$
\rho=a_u+\sum_{i>0}a_{u+i}x^i+\sum_{i>0}b_{v+i}y^i
=b_v+\sum_{i>0}a_{u+i}x^i+\sum_{i>0}b_{v+i}y^i.
$$
Since $\rho$ is a unit on a neighbourhood of $E$, possibly after 
shrinking $U$ it makes
sense to set
$$
\zeta_i=\zeta_i'/\rho,\quad i=0,1,
$$
and these local generators of $\ma{L}^{\otimes r}$ have all the 
required properties.
\end{proof}
\begin{remark}
\label{localnormal}
Set $d:=(u,v)$. Consider the surface
\[\w{S}_{u,v}=\{xs_0=\tau^{v/d} s_1,\ ys_1=\tau^{u/d} s_0\}
\subset \mathbb{A}^3_{\,x,y,\tau}\times\pr{1}_{s_0:s_1}\]
and the family of curves
$\w{f}_{u,v}\colon\w{S}_{u,v}\to\mathbb{A}^1_{\tau}$ induced by the
projection. Notice that if $d=1$, this is the normalization of
$f_{u,v}\colon S_{u,v}\to C_{u,v}$.

If the base $T$ is normal,
then $\X$ has the following stronger property:
\begin{enumerate}[(P)]
\item There exist neighbourhoods $T_0$ of $t_0$ in $T$ and
$U$ of $E$ in $\X$, and a regular map $h\colon T_0\to
\mathbb{A}^1_{\tau}$, such that $h(t_0)=0$ and $U$ is isomorphic to a
neighbourhood of the exceptional curve in the base change of the 
surface
$\w{S}_{u/d,\,v/d}$ via $h$.
\end{enumerate}
\end{remark}
\subsection{The model family of curves}
\label{fam}
Let $C$ be a nodal curve
and let $\ev^w$ be a weighted
subgraph of $\g_C$ satisfying condition (C1) of \ref{wg}.
Let $X$ be the curve obtained from $C$ by
blowing up the nodes in $\ev$.
Let $\un{\ma{C}}\to \un{D}$
be a semi-universal deformation of the curve
$C$, where $\un{D}$ is the
unit polydisc in $\C^{M}$ with coordinates
$t_1,\dotsc,t_{M}$. Observe that this deformation is universal if and
only if $C$ is stable. If this is the case, $M$ equals $3g-3$; 
otherwise it is strictly larger than
$3g-3$.

We shall construct a finite cover ${D}_{\ev^w}\to\un{D}$ and a family 
of
blow-ups $\pi\colon\X_{\ev^w}\to\ma{C}_{\ev^w}\to D_{\ev^w}$ (where
$\ma{C}_{\ev^w}=\un{\ma{C}}\times_{\un{D}} {D}_{\ev^w}$) with
central fiber $X$:
\[
{}{}\xymatrix{
{}{} {\X_{\ev^w}}\ar[r] & {\ma{C}_{\ev^w}} \ar[r]\ar[d]&
{}{}{{D}_{\ev^w}}\ar[d] \\
{}{}& {\un{\ma{C}}}\ar[r]&{\un{D}} }
\]
The family $\X_{\ev^w}\to D_{\ev^w}$
will depend on $C$ and on
the weighted graph $\ev^w$; it will be used to construct the universal
deformation of any limit $r$-th root $(X,L,\alpha)$ with weighted
graph $\ev^w$.

Denote by
$E_1, \dotsc, E_m$ the exceptional components of $X$ and by 
$n_1,\dotsc,n_m$
the corresponding nodes in $C$. For any $i=1,\dotsc,m$ let $u_i$,
$v_i$ be the weights associated to $n_i$
and let $\{t_i=0\}$ be the locus in $\un{D}$ where the node
$n_i$ persists. We write
$\un{D}=D\times D'$,
where $D$ is the unit polydisc
with coordinates $t_1,\dotsc,t_m$ and $D'$ corresponds to the 
remaining
$t_i$'s.
Consider
\[ C_{u_1,v_1}\times\cdots\times C_{u_m,v_m}
=\{\,w_{1}^{u_1}=z_{1}^{v_1},\,\dotsc\,,
w_{m}^{u_m}=z_{m}^{v_m}\,\}\subset\A^{2m}_{\,w_{1},z_{1},
\dotsc,w_{m},z_{m}}\]
and set
$${D}_{\ev^w}:=\prod_{i=1}^m C_{u_i,v_i}\times D'.$$
Consider the morphism $\prod_i C_{u_i,v_i}\to D$ defined
by
$$t_i=w_{i}z_{i}\quad\text{ for all }i=1,\dotsc,m $$
and the induced morphism
${D}_{\ev^w}\to\un{D}$.
Let
$\mathcal{C}_{\ev^w}\rightarrow {D}_{\ev^w}$ be the pull-back of the 
family
$\un{\ma{C}}\to \un{D}$ to ${D}_{\ev^w}$.

A neighbourhood $V_i$
of $n_i$ in $\ma{C}_{\ev^w}$ can be thought of as
\[V_i=\{\,xy=w_{i}z_{i},\, w_{1}^{u_1}=z_{1}^{v_1},\,\dotsc\,,
w_{m}^{u_m}=z_{m}^{v_m}\,\}
\subset\A^{M+m+2}_{\,x,y,w_{1},z_{1},
\dotsc,w_{m},z_{m},t_{m+1},\dotsc,t_{M}}.\]
Now let $\pi\colon
\mathcal{X}_{\ev^w}\rightarrow \mathcal{C}_{\ev^w}$ be the blow-up
locally described as follows:
\[U_i=\pi^{-1}(V_i)=\{\,xs_0=w_{i}s_1,\, ys_1=z_{i}s_0,\,
w_{1}^{u_1}=z_{1}^{v_1},\,\dotsc\,,
w_{m}^{u_m}=z_{m}^{v_m}\,\} \]
in $\A^{M+m+2}\times\pr{1}_{s_0:s_1}$.
Note that $\pi_{|U_i}\colon U_i\to V_i$
is the blow-up of the ideal $(x,w_{i})$,
or equivalently of $(y,z_{i})$.
Now ${\mathcal{X}}_{\ev^w}\rightarrow {D}_{\ev^w}$ is a family of
blow-ups of $\mathcal{C}_{\ev^w}\to {D}_{\ev^w}$
and has $X$ as central fiber.

Cover each
$U_i$ by the two affine open subsets $U_{i0}=\{s_1\neq
0\}$ and $U_{i1}=\{s_0\neq 0\}$. Set $\eta_0:=s_0/s_1$ and
$\eta_1:=s_1/s_0$. Then the total space $\X_{\ev^w}$ has equations
\begin{equation}
\label{eq}
w_{1}^{u_1}=z_{1}^{v_1},\,\dotsc\,,w_{m}^{u_m}=z_{m}^{v_m},\quad
\begin{cases}
x\eta_0=w_{i}\quad \text{ in $U_{i0}$, }
\cr
y\eta_1=z_{i}\quad \text{ in $U_{i1}$, }\cr
\end{cases}
\end{equation}
and $\eta_0\eta_1=1$ on $U_{i0}\cap U_{i1}$. Let $\ma{E}$ be the 
Cartier divisor given by
\begin{equation}
\label{D}
\text{$x^{u_i}$ on
$U_{i0}$ \ and \ $y^{v_i}$ on $U_{i1}$,\ \ \ for $i=1,\dotsc,m$. }
\end{equation}
On $U_{i0}\cap U_{i1}$ we have
$x^{u_i}/y^{v_i}=\eta_1^r\in\ol_{\X_{\ev^w}}^*(U_{i0}\cap U_{i1})$. 
Note that:
\begin{enumerate}[$\circ$]
\item
$\ma{E}$ is effective; its support in each $U_i$
is the inverse image via $\pi$ of the locus of
the $i$-th node $\{x=y=w_{i}=z_{i}\}$ in $V_i\subset\ma{C}_{\ev^w}$;
\item $\ma{E}_{|\w{X}}=\sum_{i=1}^m(u_ip_i+v_iq_i)$;
\item $\ol_{\X_{\ev^w}}(-\ma{E})_{|E_i}$ has degree $r$ for any 
$i=1,\dotsc,m$.
\end{enumerate}
The following result is a global version of \ref{local}.
\begin{lemma}
\label{pback}
Let $\ma{C}\to T$ be a family of nodal curves, $\ma{N}\in\Pic\ma{C}$ 
and
$(\ma{X}\to
T,\ma{L},\alpha)$ a limit $r$-th root
of $(\ma{C},\ma{N})$. Let $t_0\in T$,
and let $\ev^w$ be the weighted graph of the fiber over $t_0$.

Possibly after
shrinking $T$, there exist a fiber product
diagram
\[
{}{}\xymatrix{
{}{}{\X}\ar[r]^{\pi}\ar[d]_{k} & {\ma{C}} \ar[r]\ar[d]&T\ar[d]^{h}\\
{}{}{\X_{\ev^w}}\ar[r]&{\ma{C}_{\ev^w}}\ar[r]& {D_{\ev^w}} }
\]
where $h(t_0)=(0,\dotsc,0)$, and an isomorphism
$
\vartheta\colon \ma{L}^{\otimes r}\to \pi^*(\ma{N})\otimes 
k^*(\ol(-\ma{E}))
$
whose composition with the natural map
$
\pi^*(\ma{N})\otimes k^*(\ol(-\ma{E}))\to \pi^*(\ma{N})
$
is $\alpha$.
\end{lemma}
\begin{proof}
Call $C$ and $X$ the fibers of $\ma{C}\to T$ and $\X\to T$ over $t_0$.
Since $\ma{C}\to T$ is a deformation of $C$,
up to shrinking $T$ there is a morphism $\ell\colon T\to\un{D}$
such that $\ell(t_0)=(0,\dotsc,0)$ and
$\ma{C}\simeq\un{\ma{C}}\times_{\un{D}}T$. By
Lemma \ref{local}, $\ell$ lifts to $h\colon T\to D_{\ev^w}$
and $\ma{X}\simeq \X_{\ev^w}\times_{D_{\ev^w}}T$. We define 
$\vartheta$ to agree with $\alpha$ away from the exceptional 
components, as $k^*(\ol(-\ma{E}))$ is canonically isomorphic to 
$\ol_{\X}$ in this region. What must be seen is that $\vartheta$ can 
be extended to all of $\X$. Let $E_i$ be an exceptional component of 
$X$,
and let $\zeta_{\ma{N}}$ be a generator of $\pi^*(\ma{N})$ on a 
neighbourhood of $E_i$.
Near $E_i$, the space $\X$ is obtained from \eqref{eq} by base change 
via
$h$. Let $U_0$ and $U_1$ be the inverse images of $U_{i0}$ and 
$U_{i1}$. As in \ref{local}, we write $x$ and $y$ also to indicate 
the pull-backs of these functions to $\X$. We know from
\ref{local} that there are generators
$\zeta_0$ and $\zeta_1$ for $\ma{L}^{\otimes r}$, over $U_0$ and 
$U_1$ respectively,
with the property that
\begin{equation}
\label{compat}
\zeta_0=\eta_1^r\zeta_1\ \text{on }U_0\cap U_1,
\quad \alpha(\zeta_0)=x^u\zeta_{\ma{N}},\quad
\alpha(\zeta_1)=y^v\zeta_{\ma{N}},
\end{equation}
while $\pi^*(\ma{N})\otimes k^*(\ol(-\ma{E}))$ is generated on $U_0$ 
and
$U_1$ by $\zeta_0'=\zeta_{\ma{N}}\otimes k^*(x^u)$ and
$\zeta_1'=\zeta_{\ma{N}}\otimes k^*(y^v)$. Now we can extend 
$\vartheta$ across $E_i$
by sending $\zeta_0$ to $\zeta_0'$ and $\zeta_1$ to $\zeta_1'$; since
$\zeta_0'=\eta_1^r\zeta_1'$, the first of the identities
\eqref{compat} says that this is not ambiguous. On the other hand, 
since $\zeta_0'$ maps
to $x^u\zeta_{\ma{N}}$ in $\pi^*(\ma{N})$, and $\zeta_1'$ to 
$y^v\zeta_{\ma{N}}$, the
remaining identities say that the composition of $\vartheta$ with 
$\pi^*(\ma{N})\otimes
k^*(\ol(-\ma{E}))\to \pi^*(\ma{N})$ agrees with $\alpha$.
\end{proof}
\begin{remark}
Since $D_{\ev^w}$ equals $C_{u_1,v_1}\times\cdots\times 
C_{u_m,v_m}\times D'$, it
is normal if and only if, for every $i$, either $u_i=1$, or
$v_i=1$. This is always true if $r=2$ or $r=3$.
\end{remark}
\subsection{Construction of universal deformations}
\label{univdef}
We are now in a position to construct the universal deformation of a 
limit root.
We place ourselves in the set-up of \ref{moduli}.
Let $B$ be a scheme, $f\colon\ma{C}\to B$ a family of nodal
curves and $\ma{N}\in\Pic\ma{C}$ a line bundle of relative degree
divisible by $r$.

Fix $b_0\in B$, set $C=f^{-1}(b_0)$ and
let $\un{\ma{C}}\to\un{D}$ be a
semi-universal deformation of $C$. There exist a neighbourhood $B_0$ 
of $b_0$ in
$B$ and a morphism $B_0\to\un{D}$, $b_0\mapsto(0,\dotsc,0)$, such
that $\ma{C}_{B_0}\simeq\un{\ma{C}}\times_{\un{D}}B_0$.

Let $\xi:=(X,L,\alpha)$ be a limit $r$-th root of
$(C,N_{b_0})$ and let $\ev^w$ be the associated weighted subgraph of 
$\g_C$.
Let $\X_{\ev^w}\to D_{\ev^w}$ be the deformation of
$X$ constructed above.
Set
$$U_{\xi}:=B_0\times_{\un{D}}{D}_{\ev^w}, $$
let $u_0\in U_{\xi}$ be the point that maps to $b_0\in B_0$
and define $\ma{C}_{\xi}$,
$\ma{X}_{\xi}$ in the natural way.
We denote by $\ma{N}_{\xi}$ the pull-back of $\ma{N}$ under
$\ma{C}_{\xi}\to\ma{C}_{B_0}$;
observe that in general $\ma{N}_{\xi}$ is not the
pull-back of a line bundle on $\ma{C}_{\ev^w}$.
\[
{}{}\xymatrix{& &{\ma{N}_{\xi}}\ar[d]&&\\
{}{}{X\,}\ar@{^{(}->}[r] & {\ma{X}_{\xi}}\ar[r]^{\pi_{\xi}}\ar[d]
{}{}&{\ma{C}_{\xi}}\ar[r]\ar[d]&{U_{\xi}}\ar[r]\ar[d]&{B_0}
{}{}\ar[d]\\
{}{}&{\X_{\ev^w}}\ar[r]&{\ma{C}_{\ev^w}}\ar[r]&{{D}_{\ev^w}} 
\ar[r]&{\un{D}}}
\]
Our goal is to prove that $L$ extends to a line bundle
on $\X_{\xi}$, giving
a family of limit $r$-th roots of $\ma{N}_{\xi}$.

A key step is the next lemma, concerning automorphisms of the above
diagram. As usual, when $W\to Z$ is a morphism of schemes, we write
$\Aut_Z(W)$ to denote the group of automorphisms of $W$ over $Z$.
We shall make use of Lemma~\ref{iness}, by fixing an orientation on 
the
graph $\Sigma _X$ and hence an isomorphism
$\Aut_C(X)\simeq \ma{C}^1(\Sigma _X,\C^*)$.
\begin{lemma}
\label{G}
The action of $G = \Aut _{\un{D}}({D}_{\ev^w})$ on 
${D}_{\ev^w}$ lifts to:
\begin{enumerate}[a)]
\item
a natural action 
on
$\X_{\ev^w}$, which leaves the divisor $\ma{E}$ invariant;
\item
natural actions\ \
$G\to\Aut _{B_0}(U_{\xi}),\
G\to 
\Aut_{B_0}(\ma{C}_{\xi}),\ 
G\to\Aut_{B_0}({\ma{X}}_{\xi}),$\ \
compatible with the projections $\X_{\xi} \to\ma{C}_{\xi}\to 
U_{\xi}$.
\end{enumerate}
Moreover, the homomorphism $G\to \Aut(X)$ 
induced by
$G\to\Aut_{B_0}({\ma{X}}_{\xi})$
yields an isomorphism
$\ 
G\stackrel{\sim}{\la}\ma{C}^1(\Sigma_X,\mu_r )
\subset\Aut_C(X)$.
\end{lemma}
\begin{proof}
Automorphisms of ${D}_{\ev^w}$ over
$\un{D}$ are automorphisms
of $C_{u_1,v_1}\times\cdots\times C_{u_m,v_m}$ over $D$. Any such
automorphism has the form $(w_{i},z_{i})\mapsto(\zeta_i
w_{i},\zeta_i^{-1} z_{i})$ with $\zeta_i$ an $r$-th root of unity;
thus $G\simeq (\mu_r )^m$.

We show that $G$ acts on $\X_{\ev^w}$ and that the action on the 
central
fiber is by inessential automorphisms, giving
$G\simeq\ma{C}^1(\Sigma_X,\mu_r )\subset\Aut_C(X)$ (cf. Lemma 
\ref{iness}).
Fix $(\zeta_1,\dotsc,\zeta_m)\in(\mu_r )^m$. In each $U_{i}$
the space $\X_{\ev^w}$ has equations \eqref{eq} and
the action is given by
\begin{equation}
\label{action}\aligned
(x,\eta_0,\dots,w_{i},\dots)\mapsto
(x,\zeta_i\eta_0,\dots,\zeta_iw_{i},\dots) \
\text{ in $U_{i0}$, }\\
(y,\eta_1,\dots,z_{i},\dots)\mapsto
(y,\zeta_i^{-1}\eta_1,\dots,\zeta_i^{-1}z_{i},\dots) \
\text{ in $U_{i1}$,}
\endaligned
\end{equation}
hence the induced
automorphism of $X$ is inessential and corresponds to
$(\zeta_1,\dotsc,\zeta_m)\in\ma{C}^1(\Sigma_X,\mu_r )$. It is clear 
that $\ma{E}$ is
invariant under this action.

Now, using the
universal property of fiber products, it is immediate to see that the 
action
of $G$ lifts to $U_{\xi}$, $\ma{C}_{\xi}$, and $\X_{\xi}$ with the 
desired properties.
\end{proof}

In what follows, if $\sigma\in G$, we shall denote by 
$\sigma_{U_{\xi}}$,
$\sigma_{\ma{C}_{\xi}}$, and $\sigma_{\X_{\xi}}$ the corresponding 
automorphisms of
$U_{\xi}$, $\ma{C}_{\xi}$, and $\X_{\xi}$.
\medskip

We now address the problem of extending $L$.
Recall that there exists an effective
Cartier divisor $\ma{E}$ on $\X_{\ev^w}$ such
that $\ma{E}_{|\w{X}}=\sum_{i=1}^m(u_ip_i+v_iq_i)$ and
$\ol_{\X_{\ev^w}}(-\ma{E})_{|E_i}$ has degree $r$ for any 
$i=1,\dotsc,m$.
Denote by $\ma{G}_\xi$ the pull-back of $\ol_{\X_{\ev^w}}(-\ma{E})$ to
$\ma{X}_{\xi}$.
Then part {\it (i)} of Lemma \ref{pull} implies that there are an
inessential automorphism $\sigma$ of $X$ and an isomorphism
$
\vartheta\colon L^{\otimes
r}\to\sigma^*(\pi_{\xi}^*(\ma{N}_{\xi})\otimes\ma{G}_\xi\otimes 
\ol_{X})
$
such that $\alpha$ is the composition of $\vartheta$ with the 
pull-back of $\pi_{\xi}^*(\ma{N}_{\xi})\otimes\ma{G}_\xi\otimes 
\ol_{X}\to \pi_{\xi}^*(\ma{N}_{\xi})\otimes \ol_{X}$.
Hence, up to modifying via $\sigma$ the identification of $X$
with the central fiber of $\X_{\xi}$, we can assume that
$\alpha\colon L^{\otimes r}\to \pi_{b_0}^*(N_{b_0})$ agrees with the 
restriction to $X$ of
the natural map
$\pi_{\xi}^*(\ma{N}_{\xi}) 
\otimes\ma{G}_\xi\to\pi_{\xi}^*(\ma{N}_{\xi})$.

Shrinking $U_{\xi}$ and $B_0$, if necessary, we can extend $L$ to
$\mathcal{L}_{\xi}\in\Pic\mathcal{X}_{\xi}$ so that
$$\mathcal{L}_{\xi}^{\otimes 
r}\simeq\pi_{\xi}^*(\ma{N}_{\xi})\otimes\ma{G}_\xi\,,$$
by Remark \ref{fibrati}. Let $\alpha_{\xi}$ be the composition of 
this isomorphism with $\pi_{\xi}^*(\ma{N}_{\xi}) 
\otimes\ma{G}_\xi\to\pi_{\xi}^*(\ma{N}_{\xi})$.
Now $(\X_{\xi}\to U_{\xi},\ma{L}_{\xi},\alpha_{\xi})$ is a limit
$r$-th root of $\ma{N}_{\xi}$. Moreover, there is an isomorphism 
$\widehat{\psi}$ of limit roots between $\xi$ and the fiber of the
family $(\X_{\xi}\to U_{\xi},\ma{L}_{\xi},\alpha_{\xi})$ over $u_0\in
U_{\xi}$.
The pair $((\X_{\xi}\to 
U_{\xi},\ma{L}_{\xi},\alpha_{\xi}),\widehat{\psi})$ is a universal 
deformation for $\xi$.
\begin{prop}
\label{univ}
Let $p\colon T\to B$ be a morphism of schemes,
and $(\X\to T,\ma{L},\beta)$ a limit $r$-th root of $\ma{N}_T$.
Let $t_0\in T$ be such that $p(t_0)=b_0$, and assume that there is an 
isomorphism
$\widehat{\phi}$ of limit roots between
$\xi=(X,L,\alpha)$ and the fiber of $(\X\to T,\ma{L},\beta)$ over 
$t_0$.
Then, possibly after
shrinking $T$, the deformation $((\X\to 
T,\ma{L},\beta),\widehat{\phi})$ is isomorphic to the pull-back of 
$((\X_{\xi}\to U_{\xi},\ma{L}_{\xi},\alpha_{\xi}),\widehat{\psi})$ 
via a unique
morphism of $B$-schemes $\gamma\colon T\to U_{\xi}$ such that 
$\gamma(t_0)=u_0$; moreover, the isomorphism is unique.
\end{prop}
\begin{proof}[Proof of existence]
We know by Lemma \ref{pback} that, up to shrinking $T$, there exists
$h\colon T\to D_{\ev^w}$ such that $\ma{C}_T\simeq
\ma{C}_{\ev^w}\times_{D_{\ev^w}}T$ and
$\X\simeq\X_{\ev^w}\times_{D_{\ev^w}}T$.
Recall that $U_{\xi}=B_0\times_{\un{D}}D_{\ev^w}$, and set 
$\gamma=p\times h$. Then $\gamma\colon T\to U_{\xi}$ has the property 
that
$\gamma(t_0)=u_0$ and that there is a commutative diagram with fiber 
product squares
\begin{equation}
\label{cartsquares}
{}{}\xymatrix{
{}{}{\X}\ar[r]^{\pi}\ar[d]^{\delta} & {\ma{C}_T}\ar[r]\ar[d]^{\kappa}
{}{}& T\ar[d]^{\gamma}\ar[rd]^{p} & \\
{}{}{\X_{\xi}} \ar[r]^{\pi_{\xi}}&{\ma{C}_{\xi}}\ar[r]&
{}{} {U_{\xi}}\ar[r]&{B_0\subset B} }
\end{equation}
Let $\ma{L}'$ and $\beta'$ be the line bundle and the homomorphism
obtained by pulling back $\ma{L}_{\xi}$ and $\alpha_{\xi}$ via 
$\delta$.
Notice that $(\X,\ma{L}',\beta')$ is a limit $r$-th root of
$\ma{N}_T$, since ${\kappa}^*(\ma{N}_{\xi})=\ma{N}_T$ by the
commutativity of \eqref{cartsquares}. Moreover, it follows from Lemma 
\ref{pback} and
from the definition of $\ma{L}_\xi$ that there is an isomorphism
$\vartheta\colon \ma{L}^{\otimes r}\to{\ma{L}'}^{\otimes r}$ such that
$\beta'\circ\vartheta=\beta$.

Now let $\phi\colon X\to X_{t_0}$ and $\psi\colon X\to 
(X_{\xi})_{u_0}$ be the
isomorphisms of schemes underlying $\widehat{\phi}$ and 
$\widehat{\psi}$,
and set $L'=\phi^*(\ma{L}')$. Then 
$\sigma:=\psi^{-1}\circ\delta\circ\phi$
is an inessential automorphism of $X$ with the property that
$\sigma^*(L)=L'$, and $\vartheta$ yields an isomorphism $\chi\colon 
L^{\otimes
r}\to \sigma^*(L)^{\otimes r}$ such that
$\alpha=\sigma^*(\alpha)\circ\chi$. Then part {\it (ii)} of Lemma
\ref{pull} says that $\sigma$ belongs to $\ma{C}^1(\Sigma_X,\mu_r )$
which, by Lemma \ref{G}, can be identified with $G = \Aut
_{\un{D}}({D}_{\ev^w})$. Moreover, Lemma \ref{G}, suitably 
interpreted,
implies that $\psi\circ\sigma=\sigma_{\X_\xi}\circ\psi$. Thus, 
replacing
$\gamma$ with $\sigma_{U_\xi}^{-1}\circ\gamma$, $\kappa$ with
$\sigma_{{\ma{C}}_\xi}^{-1}\circ\kappa$, and $\delta$ with
$\sigma_{\X_\xi}^{-1}\circ\delta$, we may suppose that
$\psi=\delta\circ\phi$ (recall that, always by Lemma \ref{G}, 
$\ma{E}$ is
$G$-invariant). A consequence is that $\ma{L}$ and $\ma{L}'$ agree at
$t_0$. On the other hand, $\ma{L}^{\otimes r}$ and ${\ma{L}'}^{\otimes
r}$ are isomorphic via $\vartheta$, so Remark \ref{fibrati} implies 
that
there is an isomorphism between $\ma{L}$ and $\ma{L}'$ whose $r$-th 
tensor
power is $\vartheta$. At this point we are almost done:
$(\X,\ma{L},\beta)$ is certainly the pull-back of
$({\X}_\xi,{\ma{L}}_\xi,\alpha_\xi)$. Moreover, $\widehat{\phi}$ is 
the pull-back of
$\widehat{\psi}$, except for one little detail, namely that the
corresponding isomorphisms between $L$ and the pull-back of $\ma{L}$ 
might
not be equal. However, since the $r$-th tensor powers of these
isomorphisms are equal, the ratio between the two is an $r$-th root of
unity. To fix this it suffices to multiply the isomorphism between
$\ma{L}$ and $\ma{L}'$ by a suitable root of unity.

\bigskip
\noindent{\it Proof of uniqueness.}
Assume that the morphism $p\colon 
T\to B_0$ has another lifting $\gamma'\colon T\to U_{\xi}$ such that
$((\X,\ma{L},\beta),\widehat{\phi})$ is isomorphic to the pull-back 
of $((\X_{\xi},\ma{L}_{\xi},\alpha_{\xi}),\widehat{\psi})$ via 
$\gamma'$, and let
\begin{equation}
\label{cartsquaresbis}
{}{}\xymatrix{
{}{}{\X}\ar[r]^{\pi}\ar[d]^{\delta'} & 
{\ma{C}_T}\ar[r]\ar[d]^{\kappa'}
{}{}& T\ar[d]^{\gamma'}\ar[rd]^{p} & \\
{}{}{\X_{\xi}} \ar[r]^{\pi_{\xi}}&{\ma{C}_{\xi}}\ar[r]&
{}{} {U_{\xi}}\ar[r]&{B_0\subset B} }
\end{equation}
be the corresponding commutative diagram.
Call $h $ and $h'$ the morphisms $T\to D_{\ev^w}$ induced 
respectively by
$\gamma$ and $\gamma'$. Then
$h$ and $h'$ are both liftings of
$T\to\un{D}$, so there exists $\sigma\in G$ such that
$h'=\sigma\circ h$. As $p$ is fixed, this yields
that $\gamma'=\sigma_{U_{\xi}}\circ\gamma$, and hence also that 
$\kappa'=\sigma_{\ma{C}_{\xi}}\circ\kappa$.
If we compose the 
vertical arrows of \eqref{cartsquaresbis} with the actions of the 
inverse of $\sigma$ on the respective target spaces, we get another 
commutative diagram whose vertical arrows are $\delta'':= 
\sigma^{-1}_{\X_{\xi}}\circ\delta'$, $\kappa$, and $\gamma$. Clearly, 
$\delta''$ and $\delta$ differ at most by an automorphism of $\X$ 
over $\ma{C}_T$; in other words, there is $\iota\in 
\Aut_{\ma{C}_T}(\X)$ such that $\delta''=\delta\circ\iota$. Moreover, 
there is an isomorphism between $\ma{L}^{\otimes r}$ and 
$\iota^*(\ma{L}^{\otimes r})$ which is compatible with $\beta$. We 
need an auxiliary lemma.
\begin{lemma}
\label{autfamily}
Let $\iota\in\Aut_{\ma{C}_T}(\X)$. 
Suppose that there exists an isomorphism $\vartheta\colon 
\ma{L}^{\otimes r}\to \iota^*(\ma{L}^{\otimes r})$ such that 
$\beta=\iota^*(\beta)\circ\vartheta$. Then there exists $\tau\in G$ 
such that, possibly after shrinking $T$, 
$\delta\circ\iota=\tau_{\X_\xi}\circ\delta$.
\end{lemma}
\begin{proof}[Proof of lemma]
This is essentially a version ``with 
parameters'' of part {\it (ii)} of \ref{pull}. Let $E$ be an 
exceptional component of $X_{t_0}$. Recall from \ref{local} that 
there are open sets $U_0$ and $U_1$ covering $E$ such that $\X$ is of 
the form $x\eta_0=h_0$ in $U_0$ and of the form $x\eta_1=h_1$ in 
$U_1$, where $\eta_0\eta_1=1$ and $h_0,h_1$ are functions on a 
neighbourhood of $t_0$ such that $h_0^u=h_1^v$. Since $\iota$ is an 
automorphism over $\ma{C}_T$, $\iota^*(x)=x$ and $\iota^*(y)=y$. 
Write
$$
\iota^*(\eta_0)=\sum_{i\ge 0}a_ix^i+\sum_{i< 
0}a_i\eta_0^{-i},\quad
\iota^*(\eta_1)=\sum_{i\ge 0}b_iy^i+\sum_{i< 
0}b_i\eta_1^{-i}.
$$
From $\iota^*(x)\iota^*(\eta_0)=h_0$ we get in 
particular that $a_i=0$ for $i\ge 0$.
Similarly, $b_i=0$ for $i\ge 0$.
On the other hand, 
$a_{-1}$ and $b_{-1}$ are both different from zero, since $\iota$ 
restricts to an automorphism of $E$.
Starting from 
$\iota^*(\eta_0)\iota^*(\eta_1)=1$, another simple power series 
computation then yields that $a_i=b_i=0$ for $i<-1$, and that 
$a_{-1}b_{-1}=1$.
Thus $\iota^*(\eta_0)=k^{-1}\eta_0$ and 
$\iota^*(\eta_1)=k\eta_1$, where $k=b_{-1}$ is a unit on a 
neighbourhood of $t_0$. To prove the lemma it suffices to show that 
$k$ is an $r$-th root of unity. Now, Lemma \ref{local} says that 
there are generators $\zeta_i$ for $\ma{L}^{\otimes r}$ on $U_i$, 
$i=0,1$, and $\zeta_{\ma{N}}$ for $\pi^*(\ma{N}_T)$ on $U_0\cup U_1$, 
with the property that $\zeta_0=\eta_1^r\zeta_1$, 
$\beta(\zeta_0)=x^u\zeta_{\ma{N}}$, and 
$\beta(\zeta_1)=y^v\zeta_{\ma{N}}$. Write 
$\vartheta(\zeta_i)=f_i\iota^*(\zeta_i)$, $i=0,1$, where
$$
f_0=\sum_{i\ge 0}c_ix^i+\sum_{i< 0}c_i\eta_0^{-i},\quad
f_1=\sum_{i\ge 0}d_iy^i+\sum_{i< 0}d_i\eta_1^{-i}.
$$
From 
$\beta(\zeta_0)=x^u\zeta_{\ma{N}}$ and 
$\beta=\iota^*(\beta)\circ\vartheta$ we then get that $x^u=x^uf_0$; 
this yields that $c_i=0$ for $i>0$ and that $c_0=1$. Similarly, 
$d_i=0$ for $i>0$, and $d_0=1$. On the other hand,
$$
\eta_1^r f_1 
\iota^*(\zeta_1)=
\eta_1^r \vartheta(\zeta_1)=
\vartheta(\zeta_0)=
f_0 \iota^*(\zeta_0)=
f_0 \iota^*(\eta_1^r \zeta_1)=
k^r \eta_1^r f_0 
\iota^*(\zeta_1).
$$
This gives that 
$d_0\eta_1^r+d_{-1}\eta_1^{r-1}+\cdots= 
k^r(c_0\eta_1^r+c_{-1}\eta_1^{r+1}+\cdots)$, and hence in particular 
that $k^r=1$, as desired.
\end{proof}
We return to the proof of 
uniqueness in \ref{univ}. Lemma \ref{autfamily} implies that 
$\delta'=\sigma_{\X_\xi}\circ\tau_{\X_\xi}\circ\delta$. On the other 
hand, $\tau_{\ma{C}_\xi}\circ\kappa=\kappa$ and 
$\tau_{U_\xi}\circ\gamma=\gamma$, since $\tau_{\X_\xi}$ corresponds, 
via $\delta$, to an automorphism of $\X$ over $\ma{C}_T$. Thus, 
replacing the original $\sigma$ with $\sigma\tau$, we may in fact 
suppose that not only $\gamma'=\sigma_{U_{\xi}}\circ\gamma$ and 
$\kappa'=\sigma_{\ma{C}_{\xi}}\circ\kappa$, but also 
$\delta'=\sigma_{\X_{\xi}}\circ\delta$. In particular, 
$\sigma_{\X_{\xi}}\circ\delta\circ\phi=\psi$.
However, by the proof 
of existence, $\psi$ is
also equal to $\delta\circ\phi$, so $\sigma_{\X_{\xi}}$ acts 
trivially on the fiber of
$\X_{\xi}\to U_{\xi}$ at $u_0$. Since this action corresponds via 
$\psi$ to the action
of $\sigma$ on $X$, it follows that $\sigma=\Id$, $\gamma'=\gamma$, 
and $\delta'=\delta$. That there is a unique isomorphism 
$\ma{L}\simeq \delta^*(\ma{L}_\xi)$ which is compatible with $\beta$ 
and $\delta^*(\alpha_\xi)$ and agrees with the given one on the 
fibers above $t_0$ follows from the uniqueness part of \ref{fibrati}. 
This finishes the proof of \ref{univ}.
\end{proof}
\begin{remark}
\label{univ2}
It follows from the construction of
$(\X_{\xi}\to U_{\xi},\ma{L}_{\xi},\alpha_{\xi})$
that this family is a universal
deformation for any one of its fibers.
\end{remark}
\subsection{The moduli scheme of limit roots}
\label{proof}
We are ready to prove Theorem~\ref{coarsemoduli}.
Consider the set
\[\snbar=\coprod_{b\in B}\{\text{limit $r$-th roots of 
$(C_b,N_b)$}\}\Big/
\sim\]
where $\sim$ is the equivalence relation given by
isomorphism of limit roots.

Fix $b_0\in B$ and a limit $r$-th root $\xi=(X,L,\alpha)$ of
$(C_{b_0},N_{b_0})$. Let
$(\X_{{\xi}}\to U_{\xi},\ma{L}_{{\xi}},\alpha_{{\xi}})$ be the 
universal
deformation of $\xi$ and $B_0\subset B$ the open subset of $B$
dominated by $U_{\xi}$.

Recall that, by
Lemmas \ref{iness} and \ref{G}, there are natural identifications and 
maps:
\[\Aut(\xi)=\ma{C}^0(\Sigma_X,\mu_r)\to
\ma{C}^1(\Sigma_X,\mu_r)=G. \]
Via this homomorphism, the group 
$\Aut(\xi)$ acts on $U_\xi$, $\ma{C}_\xi$, and $\X_\xi$; if 
$\sigma\in\Aut(\xi)$, we will denote by $\sigma_{U_\xi}$, 
$\sigma_{\ma{C}_\xi}$, and $\sigma_{\X_\xi}$ the corresponding 
automorphisms of these spaces.
The universality of the family
$(\X_{{\xi}},\ma{L}_{{\xi}},\alpha_{{\xi}})$ shown in Proposition
\ref{univ} implies that, if $\sigma\in\Aut(\xi)$, then
$\sigma_{\X_{{\xi}}}^{\,*}(\ma{L}_{{\xi}})\simeq\ma{L}_{{\xi}}$ and
$\sigma_{\X_{{\xi}}}^{\,*}(\alpha_{{\xi}})\simeq\alpha_{{\xi}}$.
\begin{lemma}
\label{isoauto}
Let $b\in B_0$, let $u_1,u_2\in U_{\xi}$ be points lying over $b$, 
and let
$(X_1,L_1,\alpha_1)$, $(X_2,L_2,\alpha_2)$ be the corresponding limit
$r$-th roots of $(C_b,N_b)$.
The following are equivalent:
\begin{enumerate}[(i)]
\item
there exists $\sigma\in\Aut(\xi)$ such that
$\sigma_{U_{{\xi}}}(u_1)=u_2$;
\item there exists an isomorphism of limit roots between
$(X_1,L_1,\alpha_1)$ and $(X_2,L_2,\alpha_2)$.
\end{enumerate}
\end{lemma}
\begin{proof}
The implication $(i)\Rightarrow (ii)$ is an immediate consequence of
what we observed above:
if there exists $\sigma\in\Aut(\xi)$ such that
$\sigma_{U_{{\xi}}}(u_1)=u_2$, then
$(\sigma_{\X_{{\xi}}})_{|X_1}\colon X_1\to X_2$ is an inessential
isomorphism which induces an isomorphism of limit roots.

Conversely, assume that there is an inessential isomorphism 
$\phi\colon
X_1\to X_2$ that induces an isomorphism between the corresponding
limit roots. We are going to show first of all that $\phi=
(\sigma_{\X_{{\xi}}})_{|X_1}$ for some $\sigma\in G$, and then that
$\sigma$ belongs to the image of $\Aut(\xi)$ in $G$.
Let $d\in\un{D}$ be the image of $b$, and $d_1,d_2\in D_{\ev^w}$ the
images of $u_1,u_2$.
Recall that $G$ is the
group of automorphisms of $D_{\ev^w}$ over $\un{D}$,
and that $\un{D}=D_{\ev^w}/G$. Choose $\tau\in G$ such
that $\tau(d_1)=d_2$. Then $\tau_{U_{\xi}}(u_1)=u_2$ and
\[\rho:=\phi^{-1}\circ
(\tau_{\X_{\xi}})_{|X_1}\colon X_1\la X_1\]
is an inessential
automorphism of $X_1$. Notice that $\rho\in
\ma{C}^1(\Sigma_{X_1},\mu_r)\subset\Aut_{C_b}(X_1)$ by construction.

Set $\xi_1:=(X_1,L_1,\alpha_1)$. The universal deformation
$(\X_{{\xi_1}}\to U_{{\xi}_1},\ma{L}_{{\xi}_1},\alpha_{{\xi}_1})$
of $\xi_1$ is
naturally an open subfamily of $(\X_{\xi}\to
U_{{\xi}},\ma{L}_{{\xi}},\alpha_{{\xi}})$ (see Remark \ref{univ2}).
Applying Lemma \ref{G} with respect to $\xi_1$, we see that
there is a subgroup $G_1$ of $G$ which is
naturally identified with
$\ma{C}^1(\Sigma_{X_1},\mu_r)\subset\Aut_{C_b}(X_1)$.
If $\{j_1,\dotsc,j_k\}\subset\{1,\dotsc,m\}$ are the indices for
which the exceptional component $E_i$ of $X$ does not deform to $X_1$
(namely, in $C_b$ the $i$-th node is smoothed),
then under the isomorphism $G\simeq(\mu_r)^{m}$, $G_1$ is the
subgroup of $(\zeta_1,\dotsc,\zeta_m)$ with
$\zeta_{j_1}=\cdots=\zeta_{j_k}=1$.

Let $\tau'\in G_1$ be such that $(\tau'_{\X_{\xi}})_{|X_1}=\rho$,
and set
$\sigma:=\tau\circ(\tau')^{-1}.$
Then
$(\sigma_{\X_{\xi}})_{|X_1}=\phi$.
Notice that $\sigma_{\X_{\xi}}^{*}(\ma{G}_{\xi})\simeq\ma{G}_{\xi}$ 
and
$\sigma_{\X_{\xi}}^{*}(\pi_{{\xi}}^{*}(\ma{N}_{\xi}))\simeq
\pi_{{\xi}}^{*}(\ma{N}_{\xi})$
(this is true for all elements of $G$). In fact, the first property 
is easily deduced
from \eqref{D} and \eqref{action}; for the second property, it is
enough to observe that
$\sigma_{\ma{C}_{\xi}}^{*}(\ma{N}_{\xi})\simeq\ma{N}_{\xi}$ and that
$\pi_{\xi}\circ\sigma_{\X_{\xi}}
=\sigma_{\ma{C}_{\xi}}\circ \pi_{{\xi}}$, so
$\sigma_{\X_{\xi}}^{*}(\pi_{{\xi}}^{*}(\ma{N}_{\xi}))\simeq
\pi_{{\xi}}^{*}(\sigma_{\ma{C}_{\xi}}^{*}(\ma{N}_{\xi}))\simeq
\pi_{{\xi}}^{*}(\ma{N}_{\xi})$.

Now, since $\ma{L}_{{\xi}}^{\otimes
r}\simeq\pi_{{\xi}}^{*}(\ma{N}_{\xi})\otimes \ma{G}_{\xi}$,
we get
$\sigma_{\X_{\xi}}^{*}(\ma{L}_{{\xi}}^{\otimes
r})\simeq\ma{L}_{{\xi}}^{\otimes r}$. By construction
\[\sigma_{\X_{\xi}}^{*}(\ma{L}_{{\xi}})_{|X_1}=\phi^*(L_2)\simeq
L_1=\ma{L}_{{\xi}|X_1},\]
so we get
$\sigma_{\X_{\xi}}^{*}(\ma{L}_{{\xi}})\simeq \ma{L}_{{\xi}}$
by Remark \ref{fibrati}, and
$\sigma$ comes from $\Aut(\xi)$.
\end{proof}
\begin{proof}[Proof of Theorem \ref{coarsemoduli}]
For every $r$-th 
root $\xi$ as above we have a well-defined, injective map between 
B-sets:
\[ \beta_{\xi}\colon U_{\xi}/\Aut(\xi)\la \snbar, \]
and the image of 
$\beta_{\xi}$ inherits a complex structure from $U_{\xi}/\Aut(\xi)$. 
Since $\snbar$ is covered by subsets of the form $\im\beta_{\xi}$, in 
order for these ``charts'' to define a complex structure on $\snbar$, 
the following must hold: if $\xi_1$ and $\xi_2$ are limit roots such 
that $\im\beta_{\xi_1}$ and $\im\beta_{\xi_2}$ intersect, then 
$\beta_{\xi_2}^{-1}(\im\beta_{\xi_1})$ is open and the composition 
$\beta_{\xi_2}^{-1}\circ\beta_{\xi_1}$ is holomorphic. That this is 
true follows easily from Proposition \ref{univ}, Remark \ref{univ2}, 
and Lemma \ref{isoauto}.
In fact, choose a limit root $\eta$ corresponding to a point in the 
intersection of $\im\beta_{\xi_1}$ and $\im\beta_{\xi_2}$. Then there 
are natural (and canonically determined)
open immersions $J_i\colon U_{\eta}\hookrightarrow U_{\xi_i}$, 
compatible with
the actions of $\Aut(\eta)$ and $\Aut(\xi_i)$. Thus
$J_i$ induces a morphism $\ov{J}_i\colon 
U_{\eta}/\Aut(\eta)\hookrightarrow
U_{\xi_i}/\Aut(\xi_i)$, which is an embedding by Lemma \ref{isoauto}. 
Finally, we have
$\beta_{\eta}=\beta_{\xi_i}\circ\ov{J}_i$.

The analytic morphism $p\colon\snbar\to B$ has finite fibers. 
Proposition 
\ref{ssred} implies that $p$ is proper, hence it is a 
finite projective morphism. As $B$ is quasi-projective,
so is $\snbar$.

Finally, the scheme
$\snbar$ is a coarse moduli space for the
functor $\snf$. In fact, for any $B$-scheme $T$, the map
$\snf(T)\to\ma{H}om_B(T,\snbar)$ is naturally defined as follows. To
any limit $r$-th root $(\X\to T,\ma{L},\alpha)$ of $\ma{N}_T$, we
associate the morphism $T\to\snbar$ which is locally given by
Proposition \ref{univ}. Clearly, this morphism depends only on the
isomorphism class of $(\X\to T,\ma{L},\alpha)$.
\end{proof}
\begin{proof}[Proof of Proposition \ref{smooth}]
The morphism $p$ is smooth at a point $\xi$ if and only if
$U_{\xi}/\Aut(\xi)\simeq B_0$, that is, if and only if $\Aut(\xi)=G$. 
This
is equivalent to $b_1(\Sigma_X)=0$.
\end{proof}
\section{Geometry of the moduli space of limit roots}
\subsection{Limit roots for a fixed curve: examples}
\label{fiber}
Let $C$ be a nodal curve and $N\in\Pic C$.
We denote by
$\snc$ the zero-dimensional scheme $\ov{S}^{\,r}_{f_C}(N)$,
where $f_C\colon C\to\{pt\}$ is the trivial family.
Notice that,
by Remark \ref{remark}, for any family $f\colon\ma{C}\to B$
whose fiber over $b_0\in B$ is $C$ and for any $\ma{N}\in\Pic\ma{C}$ 
such that
$\ma{N}_{|C}=N$, the scheme-theoretical fiber over $b_0$
of the finite morphism
$\snbar\to B$ is isomorphic to $\snc$.

We have seen in Remark~\ref{graphs} that there
are $r^{b_1(\g_C)}$
weighted subgraphs of $\g_C$ satisfying conditions
(C1) and (C2).
Let $\ev^w$ be one of them. Recall that
$\ev^w$ carries the weights $u_i,v_i$ on each edge,
and that its supporting graph $\ov{\ev}$ is a
subgraph of $\g_C$. As usual, denote by $g^\nu$ the genus of the 
normalization of $C$.
The set of nodes
$\ev$ determines a partial normalization
$\w{\pi}\colon\w{X}\to C$ of $C$, and the weights determine a line
bundle $\w{\pi}^*(N)(-\sum_i(u_ip_i+v_iq_i))$ on $\w{X}$. The dual
graph of $\w{X}$ coincides with $\g_C\smallsetminus\ev$, so we have
$r^{2g^{\nu}+b_1(\g_C\smallsetminus\Delta)}$ choices
for a line bundle $\w{L}\in\Pic
\w{X}$ such that $\w{L}^{\otimes
r}\simeq\w{\pi}^*(N)(-\sum_i(u_ip_i+v_iq_i))$ ($r^{2g^{\nu}}$
choices for the pull-back of $\w{L}$ to the normalization of
$\w{X}$, and $r^{b_1(\g_C\smallsetminus\Delta)}$
choices for the gluings at the
nodes). Each choice of $\w{L}$ gives a point $(X,L,\alpha)$ in $\snc$
with weighted graph $\ev^w$.
\begin{lemma}
The geometric multiplicity of the connected component of $\snc$ 
supported on
$(X,L,\alpha)$ is
$r^{b_1(\Sigma_X)}=r^{b_1(\g_C)-b_1(\g_C\smallsetminus\Delta)}$.
\end{lemma}
\begin{proof}
Recall that $\Sigma_X$ is the graph whose
vertices are the connected components of $\w{X}$ and whose edges are
the exceptional components of $X$.
The connected component of $\snc$ corresponding to
$(X,L,\alpha)$ is isomorphic to the fiber over the origin of the 
morphism
${D}_{\ev^w}/\Aut(X,L,\alpha)\to \un{D}$ (see section 3).
If $X$ has $m$ exceptional components and $\gamma$ is
the number of connected components of $\w{X}$, then
the order
of ramification of ${D}_{\ev^w}/\Aut(X,L,\alpha)\to \un{D}$
over the origin
is $r^m/r^{\gamma-1}=r^{b_1(\Sigma_X)}$. Now notice that $\Sigma _X$
is obtained from $\g_C$ by contracting all edges
in $\g_C\smallsetminus\Delta$,
thus
$b_1(\Sigma_X)=b_1(\g_C)-b_1(\g_C\smallsetminus\Delta)$.
\end{proof}
In conclusion, using remark~\ref{graphs}, we have that
the length of $\snc$ is
\[ \sum_{\ev^w\in\partial^{-1}(\dnmodr)}
r^{2g^{\nu}+b_1(\g_C\smallsetminus\Delta)}
\cdot r^{b_1(\g_C)-b_1(\g_C\smallsetminus\Delta)}
=r^{b_1(\g_C)} \cdot
r^{2g^{\nu}+b_1(\g_C)}=r^{2g},
\]
as expected. Here, as usual, $\dnmodr$ stands for the modulo $r$ 
multidegree of $N$.
\begin{remark}
Consider a limit root $(C,L)$ having $C$ itself as underlying curve; 
then the scheme
$\snc$ is reduced at the point $(C,L)$.
\end{remark}
\begin{remark}
By construction ${D}_{\ev^w}\to \un{D}$ depends only on $C$ and on
$\ev^w$, while $\Aut(X,L,\alpha)$ depends only on $\Sigma_X$ and $r$,
by Lemma~\ref{iness}. Thus
\begin{enumerate}[$\bullet$]
\item
\emph{the scheme structure on $\snc$ depends only on $\dnmodr$.}
\end{enumerate}
\end{remark}
\begin{example}[Compact type]
Let $C=\cup_{j}C_j$ be a curve of compact type. Then
$b_1(\g_C)=0$, hence there is a unique weighted subgraph of
$\g_C$ satisfying (C1) and (C2)
(which depends on $\dnmodr$ by the previous remark).
So there is a unique blow-up $X$ of $C$ ($X$ is of
compact type)
with a divisor $D=\sum_i(u_ip_i+v_iq_i)$, such that every limit
root in $\snc$ is of type $(X,L,\alpha)$
with $L_{|\w{X}}^{\otimes
r}\simeq\w{\pi}^*(N)(-D)$.
In particular, $X=C$ if and only if
$\dnmodr=0$.
Notice also that $\snc$ is reduced.
\end{example}
\begin{example}
\label{es1}
Let $C$ be a curve of genus $g$ with two smooth components $C_1$, 
$C_2$
and three nodes $n_1,n_2,n_3$.
\begin{center}
{}{}\scalebox{0.40}{\includegraphics{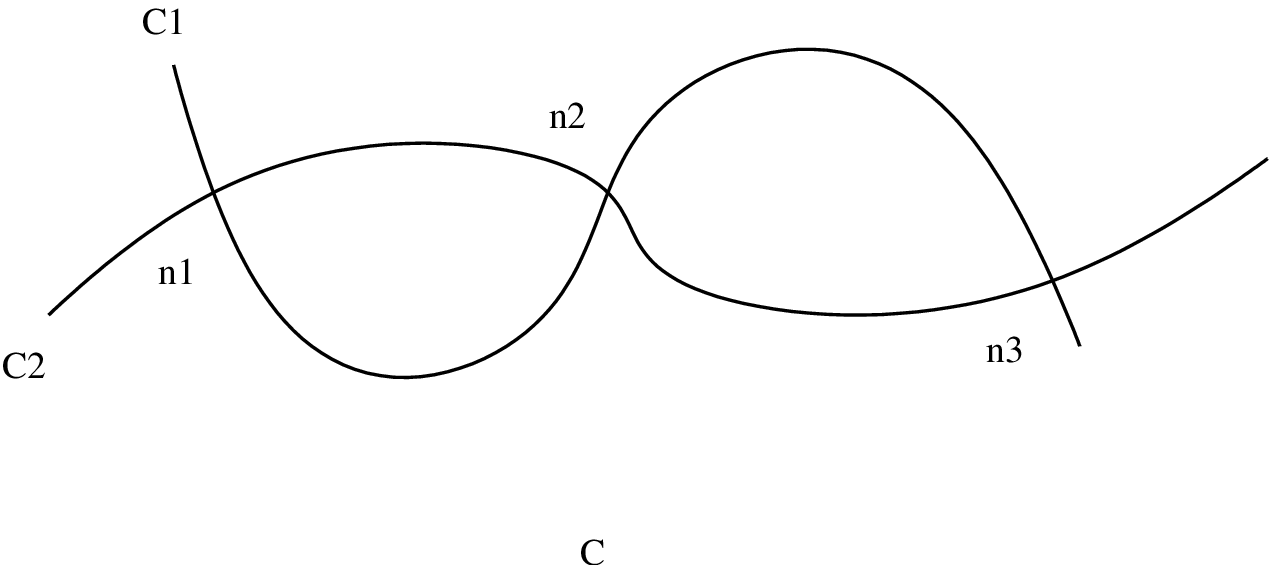}}\hspace{80pt}
{}{}\scalebox{0.35}{\includegraphics{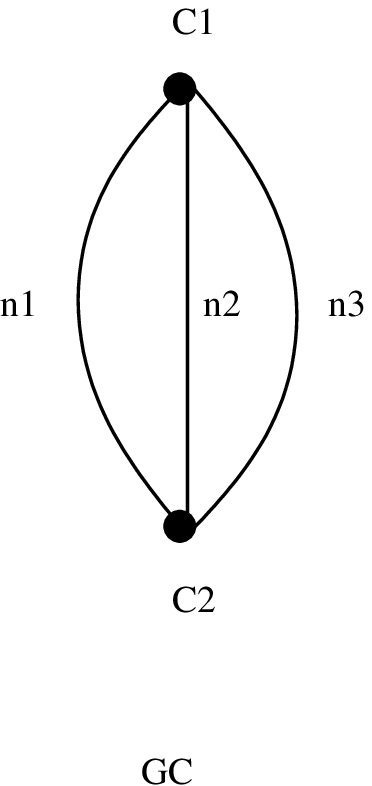}}
\end{center}
We describe $\ov{S}^3_C(\ol_C)$.
Since $b_1(\g_C)=2$, $\g_C$ must have $3^2=9$ weighted
subgraphs satisfying (C1) and (C2).
One is $\emptyset$. This
choice corresponds to points of type $(C,L)$ with $L^{\otimes
3}\simeq\ma{O}_C$; there are $9\cdot 3^{2g^{\nu}}$ of them
(with ${g^{\nu}}=g-2$), all reduced points
in $\ov{S}^3_C(\ol_C)$.

For each $\Delta$ consisting of two edges $\{n_j,n_k\}$
(where $j,k\in\{1,2,3\}$), we have two possible choices for the
weights:
\begin{center}
{}{}\scalebox{0.35}{\includegraphics{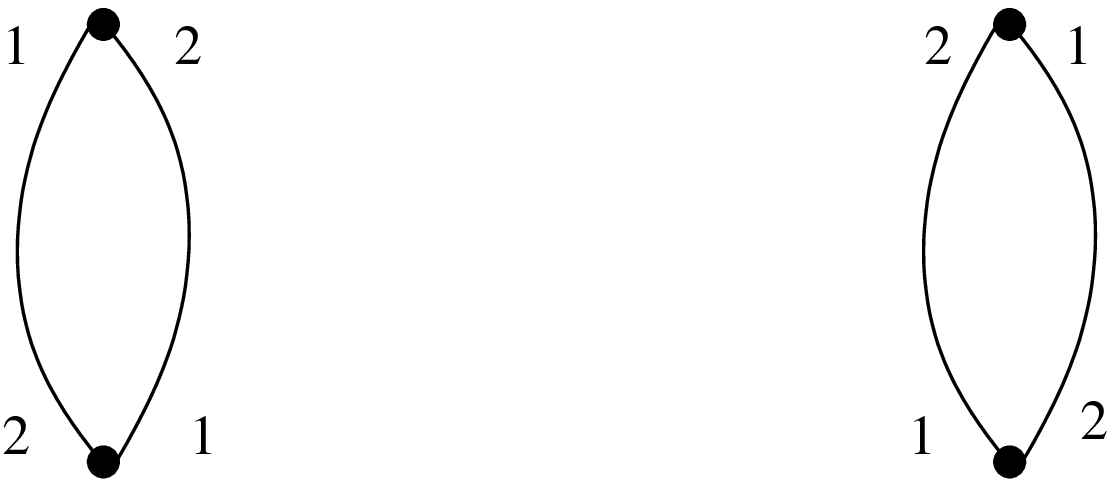}}
\end{center}
Denote by $X_{j,k}$ the curve
obtained blowing up $\ev$.
\begin{center}
{}{}\scalebox{0.40}{\includegraphics{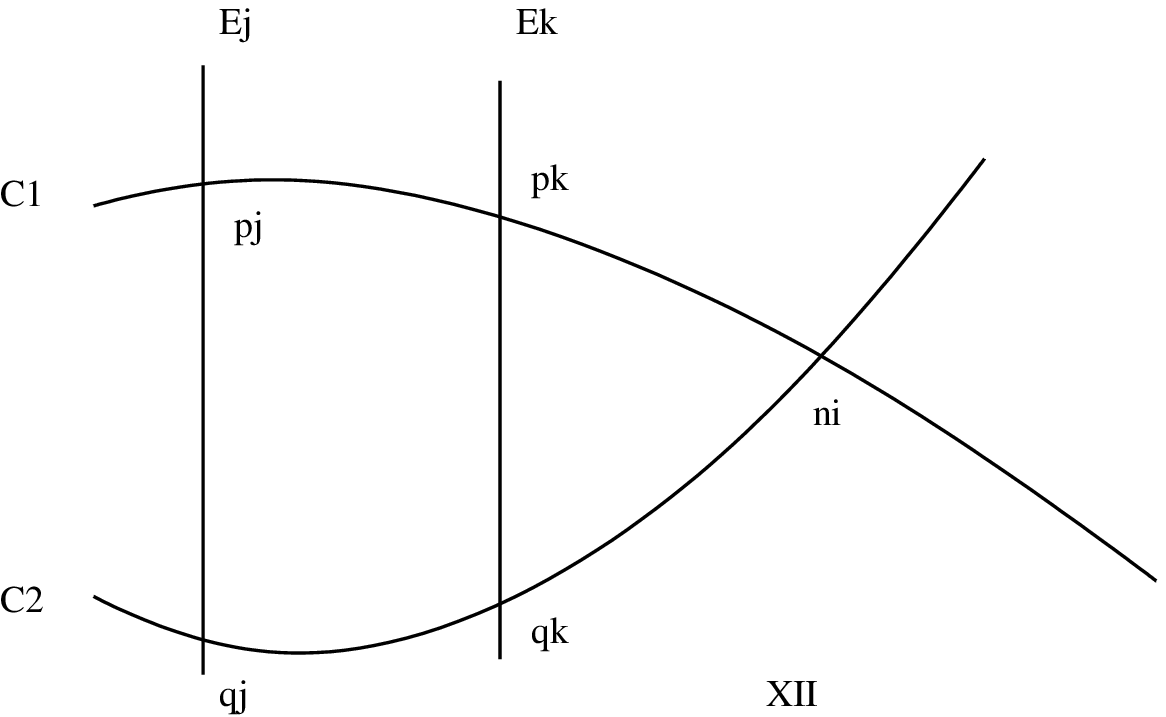}}
\end{center}
We get points of type $(X_{j,k},L)$ where
$L_{|\w{X}_{j,k}}^{\otimes
3}\simeq\ma{O}_{\w{X}_{j,k}}(-2p_j-q_j-p_k-2q_k)$
or $L_{|\w{X}_{j,k}}^{\otimes
3}\simeq\ma{O}_{\w{X}_{j,k}}(-p_j-2q_j-2p_k-q_k)$, depending on
the weights. Since $\g_C\smallsetminus\Delta$ is contractible, we
have $b_1(\g_C\smallsetminus\Delta)=0$, so these points have 
multiplicity 9 in
$\ov{S}^3_C(\ol_C)$; there are $6\cdot 3^{2g^{\nu}}$ of them.

Finally, when $\Delta=\{n_1,n_2,n_3\}$,
we have two possible choices for the weights:
\begin{center}
{}{}\scalebox{0.35}{\includegraphics{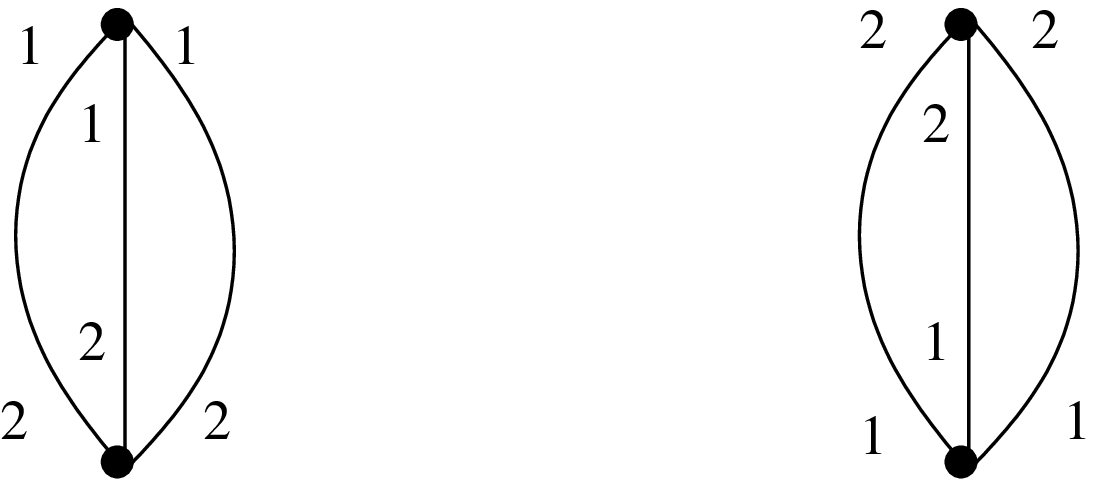}}
\end{center}
Denote by $X_{1,2,3}$ the corresponding curve.
\begin{center}
{}{}\scalebox{0.40}{\includegraphics{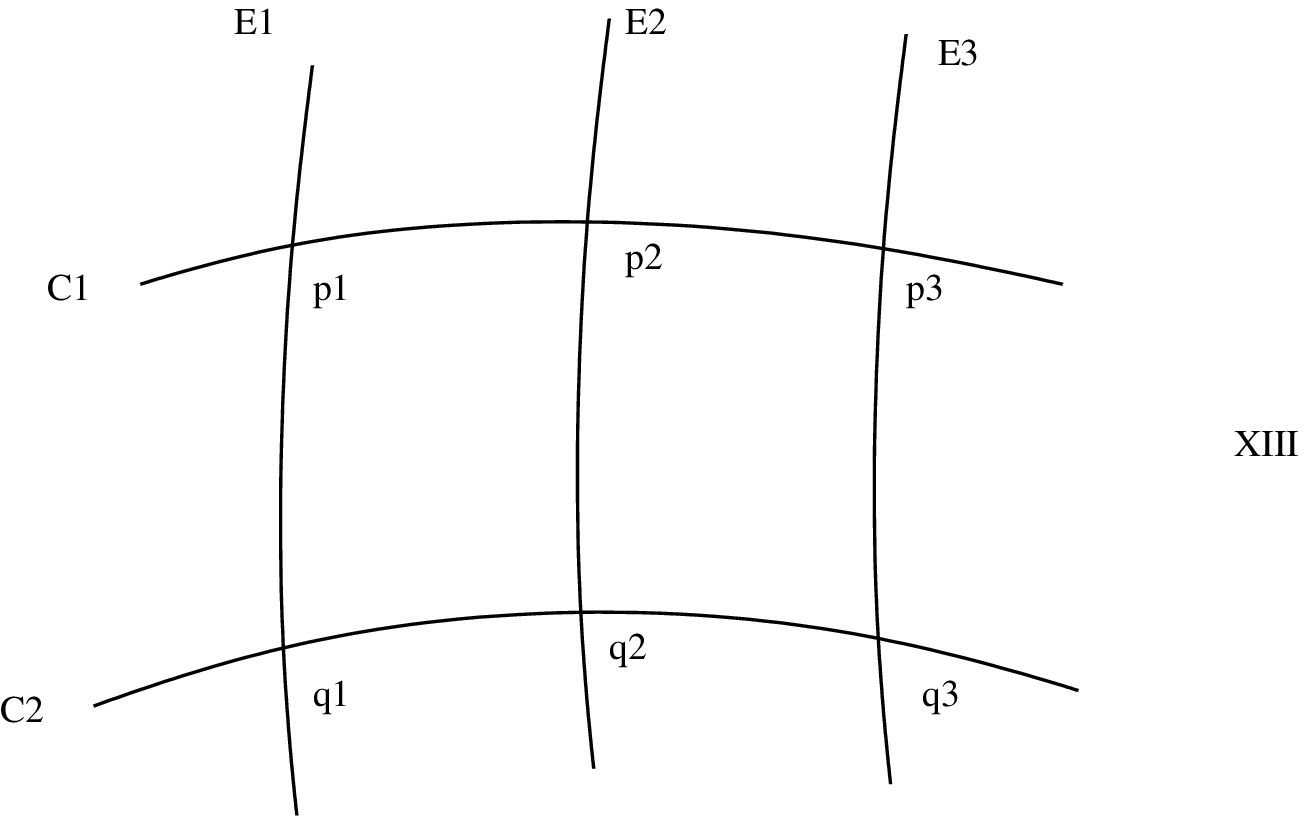}}
\end{center}
We get points $(X_{1,2,3},L)$ where $L_{|\w{X}_{1,2,3}}^{\otimes
3}\simeq\ma{O}_{\w{X}_{1,2,3}}(-p_1-p_2-p_3-2q_1-2q_2-2q_3)$
or $L_{|\w{X}_{1,2,3}}^{\otimes
3}\simeq\ma{O}_{\w{X}_{1,2,3}}(-2p_1-2p_2-2p_3-q_1-q_2-q_3)$. There 
are $2\cdot 3^{2g^{\nu}}$ of these points, each of them with 
multiplicity 9 in $\ov{S}^3_C(\ol_C)$.
\end{example}

\subsection{Higher spin curves}
\label{spin}
Let $\mgnbar$ be the moduli space, and $\ov{\ma{M}}_{g,n}$ the 
moduli stack, of stable $n$-pointed curves of genus $g$. Let 
$\omega_{u_n}$ be the relative dualizing sheaf, and 
$\Sigma_1,\dots,\Sigma_n$ the tautological sections, of the universal 
curve $u_n\colon \ov{\ma{C}}_{g,n}\to\ov{\ma{M}}_{g,n}$.

Fix integers $l,m_1,\dotsc,m_n$ such that
$r\,|\,l(2g-2)+m_1+\cdots+m_n$, and set $\un{m}:=\{m_1,\dotsc,m_n\}$
and
\[\omega^l_{\un{m}}:=\omega^{\otimes l}_{u_n}
\,(\,m_1\Sigma_1+\cdots 
+m_n\Sigma_n\,)\
\in\ \Pic\ov{\ma{C}}_{g,n}\,.\]
We introduce
a contravariant functor
$$\sgnf \colon \{\text {schemes}\}\la
\{{\rm {sets}}\},$$
encoding the moduli problem for
limit $r$-th roots of $\omega^l_{\un{m}}$ over $\ov{\ma{M}}_{g,n}$, 
as follows. Look at pairs $(f\colon\ma{C}\to B,\xi)$ consisting of a 
family $f\colon \ma{C}\to B$ of $n$-pointed genus $g$ curves and a 
limit root $\xi$ of
$(\ma{C},\omega_f^{\otimes l}(\sum m_i\sigma_i))$, where
$\sigma_1,\dots,\sigma_n$ are the tautological sections of $f$ and
$\omega_f$ is the relative dualizing sheaf. Two such pairs $(f\colon 
\ma{C}\to B,\xi)$ and
$(f'\colon \ma{C}'\to B,\xi')$ are considered equivalent if there is 
an
isomorphism $\kappa\colon \ma{C}\to \ma{C}'$ of families of 
$n$-pointed curves
and an isomorphism, as limit roots of
$(\ma{C},\omega_f^{\otimes l}(\sum m_i\sigma_i))$, between
$\xi$ and $\kappa^*(\xi')$ (observe that, if 
$\sigma'_1,\dots,\sigma'_n$
are the tautological sections of $f'$, then $\omega_f^{\otimes
l}(\sum m_i\sigma_i)$ and $\kappa^*(\omega_{f'}^{\otimes
l}(\sum m_i\sigma'_i))$ are canonically isomorphic). Notice that this
is clearly a coarser equivalence relation than the one given by 
isomorphism
of limit roots over a fixed family of curves.

Then, for every scheme $B$, we define $\sgnf(B)$ to be the set of all
equivalence classes of pairs $(f\colon \ma{C}\to B,\xi)$ as above. 
Arguing
exactly as in the proof of Theorem
\ref{coarsemoduli} we obtain:
\begin{thm}
The functor $\sgnf$ is coarsely represented by a projective scheme
$\sgnbar$, finite over $\mgnbar$.
\end{thm}
\noindent We call $\ov{\ph}\colon\sgnbar\to\mgnbar$ the structure 
morphism.
For any stable $n$-pointed curve $(C,\sigma_1,\dotsc,\sigma_n)$,
the fiber of $\ov{\ph}$
over $[C]$ is
\[\ov{S}^{\,r}_C(\omega_C^{\otimes l}(m_1\sigma_1+\cdots
+m_n\sigma_n))\,/
\Aut(C,\sigma_1,\dotsc,\sigma_n).\]

When we consider stable curves without markings (i.e., when $n=0$), 
we denote by
$\sgbar$ the corresponding moduli space, and by $\sbarC$ the fiber of
$\ov{\ph}$ over $[C]\in\mgbar$.
The scheme $\ov{S}^{\,2,1}_g$ is the moduli space of spin curves,
constructed in \cite{cornalba1}.

Denote by $\sgn$ the open subscheme of
$\sgnbar$ parametrizing points $(C,L)$ with
$C$ stable, and let $\ph\colon \sgn\to\mgnbar$ be the
natural morphism.

Pairs $(C,L)$ with $C$ a smooth curve and $L\in\Pic C$ such that
$L^{\otimes r}\simeq\omega_C$ are generally called \emph{$r$-spin
curves}.  Their moduli space admits a natural stratification by loci parametrising
line bundles $L$ with increasing $h^0(C,L)$, so that the largest locus
corresponds to ``effective $r$-spin curves"; we refer to the recent paper of A. Polishchuk  \cite{poli}
for  more
details and  open problems related to this interesting aspect.

Several modular compactifications of the moduli space of $r$-spin
curves, and more generally of $\sgn$,
have been introduced by T.\ J.\ Jarvis \cite{jarvis1,jarvis2,jarvis3},
by means of rank 1, torsion-free sheaves.

The next statement summarizes well-known results relating line 
bundles and torsion-free sheaves
of rank 1.
\begin{prop}
\label{sheaves}
Let $B$ be an integral scheme and $f\colon\ma{C}\to B$ a family of 
nodal
curves.
\begin{enumerate}[(I)]
\item
Let
$\pi\colon\ma{X}\to\ma{C}$ be a family of blow-ups of $\ma{C}$
and let $\ma{L}\in\Pic\ma{X}$ be a line bundle having degree 1 on 
every
exceptional component. Then $\pi_*(\ma{L})$ is a relatively 
torsion-free sheaf of rank 1, flat over $B$.
\item
Conversely, suppose that $\ma{F}$ is a relatively torsion-free sheaf 
of
rank 1 on $\ma{C}$, flat over $B$. Then there exist a family
$\pi\colon\X\to\ma{C}$ of blow-ups of $\ma{C}$
and a line bundle $\ma{L}\in\Pic\X$
having degree 1 on all exceptional components, such that
$\ma{F}\simeq\pi_*(\ma{L})$.
\item
Let $\pi\colon\ma{X}\to\ma{C}$, $\pi'\colon\ma{X}'\to\ma{C}$ be
families of blow-ups of
$\ma{C}$ and $\ma{L}\in\Pic\ma{X}$,
$\ma{L}'\in\Pic\ma{X}'$ line bundles having degree 1 on every
exceptional component. Then $\pi_*(\ma{L})\simeq\pi'_*(\ma{L}')$
if and only if there exists an
isomorphism $\sigma\colon\X\stackrel{\sim}{\to}
\X'$ over $\ma{C}$ such that
$\ma{L}\simeq\sigma^*(\ma{L}')$.
\end{enumerate}
\end{prop}
\begin{proof}
{\em (I)} is a local property, see \cite{jarvis1}, Proposition
3.1.2.

{\em (II)} is a consequence of G.\ Faltings' local characterization of
torsion-free sheaves \cite{faltings};
see \cite{jarvis1}, \S 3.1 and 3.2, in
particular Theorem 3.2.2.

{\em (III)} The ``if'' part is clear. To prove the converse, it 
suffices to notice that there is an isomorphism $\tau\colon \X\to 
\operatorname{Proj}(\bigoplus_{k\ge 
0}\operatorname{Sym}^k(\pi_*(\ma{L})))$ of spaces over $\ma{C}$, with 
the property that $\tau^*(\ma{O}(1))$ is isomorphic to $\ma{L}$, and 
similarly for $\X'$ and $\ma{L}'$.
\end{proof}
\noindent
Jarvis introduces
the notion of $r$-th root sheaf of a line bundle
$N\in\Pic C$, which is a rank 1,
torsion-free sheaf $\mathcal{F}$ on $C$ together with a homomorphism
$\beta\colon\mathcal{F}^{\otimes r}\to N$ with the following
properties (see \cite{jarvis3}, Definition
2.3):
\begin{enumerate}[$\circ$]
\item $r\deg\mathcal{F}=\deg N$;
\item $\beta$
is an isomorphism on the locus where $\ma{F}$ is locally free;
\item at any node of $C$ where $\ma{F}$ is not locally free, the
length of the cokernel of $\beta$ is $r-1$.
\end{enumerate}
Such an $r$-th root sheaf is always of the form $\pi_*(L)$, where
$X$ is the curve obtained from $C$ by
blowing-up the nodes
where $\ma{F}$ is not locally free, $\pi\colon X\to C$ is the
natural morphism and $(X,L,\ov{\beta})$ is a limit $r$-th root of
$(C,N)$.
It is easy to see, using Proposition \ref{sheaves}, that
the moduli functor of $r$-th root sheaves of
$\omega^l_{\un{m}}$ over $\mgnbar$
is isomorphic to
$\sgnf$
(see \cite{jarvis1}, 3.1 and 4.2.2, and \cite{jarvis3}, 2.2.2 and page
38). Hence we have the following:
\begin{thm}
\label{Jarvis}
The scheme
$\sgnbar$ is isomorphic over $\mgnbar$ to the coarse moduli
space of the stack
${\text{\scshape Root}}_{g,n}^{1/r}(\omega^l_{\un{m}})$
(in the notation of \cite{jarvis2})
of $r$-th root sheaves of $\omega^l_{\un{m}}$.
\end{thm}
Observe that, exactly as in \cite{jarvis2}, one could construct 
different
(less singular) compactifications of $\sgn$ (and more
generally of $\sn$), in two ways:
\begin{enumerate}[(a)]
\item requiring, in the definition of a family of limit roots,
that property (P) of Remark \ref{localnormal} holds; this amounts to
considering only ``not too singular'' families of limit roots;
\item attaching to a limit $r$-th
root the data of a set of limit $d$-th roots for any positive $d$
dividing $r$, i.\ e.\ considering ``coherent nets of
roots''. Set-theoretically, this
amounts to attaching to a limit
$r$-th root $\xi=(X,L,\alpha)$ the data of the
gluings for
$\xi^{m}$ for any integer $m$ dividing $r$.
\end{enumerate}
\section{Higher spin curves in the universal Picard scheme}
\subsection{Compactifying the Picard functor}
Assume $g\geq 3$ and
let $P_{\,d,g}\to\mg^0$
be the degree $d$ universal Picard variety\footnote{$\mg^0$ (and 
likewise $\mgbar^0$) stands for the locus of curves with trivial 
automorphism group.},
whose fiber over a smooth curve $C$ is $\Pic ^dC$, the variety of 
line bundles
of degree $d$ on $C$. Let
$\pdg\to\mgbar$ be the modular compactification constructed
in \cite{caporaso};
$\pdg$ is an integral, normal projective scheme.
Its fiber over $C\in \mgbar ^0$ is denoted by $\PC$ and gives a
compactification of (a finite number of copies of)
the generalized Jacobian of $C$.
The moduli properties of $\pdg$
are described in
\cite{caporaso}, section 8.
Its boundary points are of the same nature as limit roots: they 
correspond
to line bundles on quasistable
curves\footnote{A \emph{quasistable curve} $X$ is a blow-up of a 
stable curve $C$.
We call $C$ the \emph{stable model} of $X$.}, having degree $1$ on
exceptional
components.

The moduli space of higher spin curves over $\mg^0$ naturally embeds 
in $P_{\,d,g}$
(see \ref{new});
its projective closure in $\pdg$ is thus a new compactification which
we shall now study.

We recall from \cite{caporaso} a few basic properties of $\pdg$,
which, as $\mgbar$, is constructed by means of Geometric Invariant 
Theory
(see \cite{gieseker}; in this section
we freely use the language of \cite{GIT}).
We begin with some useful conventions and numerical preliminaries.
Let $d$ be an integer;
we denote by $\un{d}$ an element of $\Z ^{\gamma}$ whose entries
add up to $d$, that is
$\un{d} =\{d_1,\dotsc, d_{\gamma}\}$ and $\sum d_i =d$.
We say that $\un{d}$ is divisible by some integer $r$ if all of
its entries are;
if that is the case, we write $\un{d}\equiv 0\pmod r$.
If $q$ is a rational number, we set
$q\un{d} =\{qd_1,\dotsc, qd_{\gamma}\}$ and we say that $q\un{d}$
is integer if $qd_i$ is integer
for every $i$.

For any nodal curve $X$, we denote by $X_1,\dotsc, X_{\gamma}$ its
irreducible components
and set $k_i=\#X_i\cap \overline{X \smallsetminus X_i}$.
Let $M$ be a line bundle of degree $d$ on $X$; denote by 
$\un{d}=\un{\deg}\ M =\{d_1,\dotsc, d_{\gamma}\}$ its multidegree, 
where $d_i :=\deg_{X_i}M$.
For the dualizing sheaf of $X$ we use the special notation
$w_i:=\deg _{X_i} \omega _X = 2g(X_i)-2+k_i$
and
$$
\un{w}_X:=\un{\deg}\ \omega _X=\{w_1,\dotsc, w_{\gamma}\}.
$$
Similarly, for a subcurve $Z$ of $X$, we set $k_Z:=\#Z\cap \overline{X
\smallsetminus Z}$,
$w_Z:=\deg _Z\omega _X= 2g(Z)-2+k_Z$ and $d_Z :=\deg _Z M$.

Let $X$ be a quasistable curve.
As usual, we denote by $E_i$ its exceptional
components and set
$\w{X}=\ov{X\smallsetminus\cup_iE_i}$ (thus
$\w{X}$ is a partial normalization
of $C$).
\begin{defi}
\label{defBI}
Let $X$ be a quasistable curve and $M\in\Pic X$.
We shall say that the multidegree $\un{d}$ of $M$ is {\em {balanced}}
if $\deg_E M=1$ on every exceptional component $E$ of $X$
and if,
for every subcurve $Z$ of $X$, the following inequality
(known as Basic Inequality) holds:
$$
\left|d_Z-d\frac{w_Z}{2g-2}\right|\leq \frac{k_Z}{2}.
$$
We shall say that $\un{d}$ is {\em{stably balanced}} if it is
balanced and if for every subcurve $Z$ of $X$ such that
$$
d_Z-d\frac{w_Z}{2g-2}=- \frac{k_Z}{2},
$$
we have that $Z$ contains $\w{X}$.

Similarly, we shall say that the line bundle $M$ is {\em {balanced}}
(or {\em {stably balanced}})
if its multidegree is.
\end{defi}
It is easy to see that if $\un{d}$ is balanced and $Z$ is a connected
component of $\w{X}$, then $d_Z-d\frac{w_Z}{2g-2}=- \frac{k_Z}{2}$.
Hence if $\un{d}$ is stably balanced, $\w{X}$ is connected.

We shall need the following elementary
\begin{lemma}
\label{bal}Let $X$ be a quasistable curve.
\begin {enumerate}[{\it (i)}]
\item
\label{baldual}
Let $M\in \Pic X$ be such that
$M^{\otimes r}\simeq \omega _X^{\otimes l}$ for some integers $r$ and
$l$; then $M$ is stably balanced.
\item
\label{baltrans}
$M\in \Pic X$ is balanced
(stably balanced)
if and only if $M\otimes \omega _X^{\otimes
t}$
is balanced (stably balanced) for some integer $t$.
\item
\label{balirr}
Assume that the stable model of $X$ is irreducible and let $M\in\Pic
X$ be a line bundle having degree 1 on all exceptional
components. Then $M$ is stably balanced.
\end{enumerate}
\end{lemma}
The goal of the following definitions is
to measure the non separatedness of
the Picard functor.
\begin{defi}
\label{twister}
Let $X$ be a nodal curve and
$T\in \Pic X$; we say that $T$ is a {\em twister} if there exists
a one parameter smoothing
$\ma{X}\to S$ of $X$ such that
$$
T\simeq {\mathcal O}_{\ma{X}}
(\sum_{i=1}^{\gamma} a_iX_i)\otimes {\mathcal O}_{X}
$$
for some Cartier divisor $\sum_{i=1}^{\gamma} a_iX_i$ on $\ma{X}$
(where $a_i \in \Z$).
The set of all twisters of $X$ is denoted by $\Tw (X)$.
\end{defi}
\begin{defi}
\label{GITeq}
Consider two balanced line bundles $M\in\Pic X$ and $M'\in\Pic X'$,
where $X$ and $X'$ are quasistable curves.

We say that
$M$ and $M'$
are {\em {equivalent}} if there exists a semistable curve $Y$
dominating\footnote{We say that a semistable curve $Y$ dominates a 
nodal curve $X$ if it can be obtained from $X$ by a finite sequence 
of blow-ups.} both $X$ and $X'$ and a twister
$T\in \Tw (Y)$ such that, denoting by $M_Y$ and $M_Y'$
the pull-backs of $M$ and $M'$ to $Y$, we have
$$
M_Y'\simeq M_Y\otimes T.
$$
In particular, for any automorphism $\sigma$ of $X$, $\sigma^*(M)$
is equivalent to $M$.
\end{defi}
To say that $M$ and $M'$ are equivalent is to say that their 
pull-backs to $Y$
are both limits of the same family of line bundles.
More precisely, let ${\mathcal Y} \la S$ be a one-parameter smoothing 
of $Y$
such that the twister $T$ of Definition~\ref{GITeq} is
$T= \ol _{\mathcal Y}(D)\otimes \ol _Y$, with $\Supp D \subset Y$. Let
${\mathcal M}\in \Pic {\mathcal Y}$ be such that ${\mathcal 
M}_{|Y}\simeq M_Y$; then
${\mathcal M}\otimes \ol _{\mathcal Y}(D) \otimes \ol _Y\simeq M_Y'$ 
so that
$M_Y$ and $M_Y'$ are both limits of the family corresponding to 
${\mathcal M}^*$
(recall that the ``$*$" denotes restriction away
from the central fiber).
Therefore, the line bundles
$M_Y$ and $M_Y'$ (and likewise $M$ and $M'$) must be identified in
any proper completion of the Picard scheme over $\mgbar$.
\begin{remark}
\label{reseq}
Let $M,M'\in\Pic X$
be stably balanced. Then a simple numerical checking shows that they 
are
equivalent if and only if there exists an automorphism $\sigma$ of 
$X$ such that
$M'\simeq \sigma^*(M)$.
If the stable model of $X$ has no non-trivial automorphisms, then $M$
and $M'$ are equivalent if and only if $M_{|\w{X}}\simeq M'_{|\w{X}}$.
\end{remark}
Fix a large $d$
($d\geq 20(g-1)$ would work), set $s=d-g$ and consider the Hilbert
scheme $\Hilb$
of connected curves of degree $d$ and genus $g$ in $\pr{s}$; the group
$SL(s+1)$
naturally acts on $\Hilb$ .
Recall that there exist linearizations for such
an action such that the following facts hold.
\begin{thm}[\cite{caporaso}]
\label{BI}
Let $X\subset\pr{s}$ be a
connected curve of genus $g$.
\begin{enumerate}[(1)]
\item
The Hilbert point of $X$
is GIT-semistable if and only if $X$ is quasistable and
$\mathcal{O}_X(1)$ is balanced.
\item
The Hilbert point of $X$ is GIT-stable if and only if $X$ is
quasistable
and $\mathcal{O}_X(1)$ is stably balanced.
\item
Assume that the Hilbert points of $X$ and $X'\subset \pr{s}$ are
GIT-semistable. Then they are GIT-equivalent if and only if
$\ol_X(1)$ and
$\ol_{X'}(1)$ are equivalent.
\end{enumerate}
\end{thm}
\begin{proof}
We have assembled together,
for convenience, various results from \cite{caporaso}:
for {\it (1)} see Propositions 3.1 and 6.1 (the ``only if'' part is 
due to D.\ Gieseker \cite{gieseker}); for {\it (2)} see Lemma 6.1.
The ``only if" of part {\it (3)} is Lemma 5.2;
the ``if" is a direct consequence of the existence
of the GIT quotient as a projective scheme (see below).
\end{proof}
Let then $\pdg$ be the
GIT-quotient of the locus of all GIT-semistable
points in $\Hilb$;
it is a normal, integral scheme, flat over $\mgbar ^0$
(Theorem 6.1 in \cite{caporaso}).

The above result permits a modular interpretation of $\pdg$ as a 
modular
compactification of the universal Picard variety.
Let $X$ be a quasistable curve and let $M\in\Pic X$ be very ample of 
degree $d$.
By Theorem~\ref{BI},
if we embed $X$ in $\pr{s}$ by $M$,
the Hilbert point of $X$ is GIT-semistable if and only if $M$ is
balanced.
Thus $\pdg$ parametrizes equivalence classes of
balanced line bundles on quasistable curves.

Associating to $M\in\Pic X$ the stable model $C$ of $X$ yields the
natural morphism $\pdg\to\mgbar$; $\PC$ denotes its fiber over
$[C]\in\mgbar$.

Recall (Theorem 6.1 in \cite{caporaso}) that $\PC$ is a finite union 
of
$g$-dimensional irreducible
components, one for each
stably balanced multidegree (of degree $d$)
on $C$:
$$
\PC =\bigcup_{{\un{d}\ \text{\tiny stably}\atop \text{\tiny 
balanced}}} \PCd.
$$
Moreover, if
$C$ has no non-trivial automorphisms, every component of $\PC$
contains a copy of the generalized Jacobian
of $C$ as a dense open subscheme:
$$
\Pic ^{\, \un{d}}C \subset \PCd,
$$
where $\Pic ^{\,\un{d}}C$ is the variety of line bundles on $C$ having
multidegree $\un{d}$.
If $C$ is irreducible or of compact type, then $\PC$ is irreducible.
\begin{remark}
\label{highdeg}
By \cite{caporaso},
Lemma 8.1, for any integer $t$, $\pdg$
and
$\overline{P}_{d+t(2g-2),\,g}$ are isomorphic over $\mgbar$
via the
morphism
$M\mapsto M\otimes \omega_X^{\otimes t}$ for any $M\in\Pic X$. Hence 
we
can assume
that $\pdg$ exists for any $d\in\Z$, even if it is
constructed only
for $d\gg 0$. Moreover, $M$ is balanced, or stably balanced,
if and only if $M\otimes \omega_X^{\otimes
t}$ is (by part {\it (\ref{baltrans})} of Lemma \ref{bal}).
\end{remark}
\subsection{A new compactification of higher spin curves}
Let $C$ be a stable curve
and
$L\in\Pic C$ an $r$-th root of
$\omega _C^{\otimes l}$;
then,
by part {\it (\ref{baldual})} of Lemma~\ref{bal}, $L$ is stably
balanced, therefore, by part {\it (2)} of \ref{BI}, it is identified 
with
a point of $\ov{P}_{l(2g-2)/r,\,g}$.
More precisely, consistently with section~\ref{spin},
we denote by $\sg$
the open subscheme of
$\sgbar$ corresponding to pairs $(C,L)$ as above.
Obviously, $\sg$ is the coarse moduli space for the
similarly defined functor
${\mathcal S}^{\,r,\,l}_g$.
We have:
\begin{lemmadef}
\label{new}
There is a natural
embedding
$$
\chi\colon\sg\hookrightarrow
\ov{P}_{l(2g-2)/r,\,g}.
$$
We define $\sgtil$ to be the closure,
in the projective variety
$\ov{P}_{l(2g-2)/r,\,g}$, of the image of $\chi$, and
$$
\widehat{\ph}\colon\sgtil\la\mgbar
$$
to be the natural morphism.
\end{lemmadef}
\begin{proof}
The set-theoretic inclusion has already been defined before the
statement. To complete the proof, notice that there is an obvious
inclusion of functors
$$
{{\mathcal S}^{\,r,\,l}_g}\hookrightarrow
\ov{{\mathcal P}}_{l(2g-2)/r,\,g}
$$
and hence the desired embedding follows from \cite{caporaso}, 
Proposition
8.1 (1).
\end{proof}
\noindent This
defines a new compactification of the moduli space of
higher spin curves, which we shall compare with the previously
constructed $\sgbar$. As for
$\sgbar$, the
boundary points of $\sgtil$ correspond to line bundles
on quasistable curves, having degree 1 on exceptional components.
Here is a more precise description.
\begin{prop}
\label{stil}
The points of
$\sgtil$ are in
bijective correspondence with equivalence classes of
balanced line bundles $M\in\Pic X$ such that
$X$ is a quasistable curve of genus $g$
and there exists a twister $T$ on $X$ for which the following
relation holds:
$$M^{\otimes r}\simeq\omega_X^{\otimes l}\otimes T.$$
\end{prop}
\begin{proof}
Let $\lambda\in \sgtil$.
Since $\lambda$ belongs to the closure of (smooth) $r$-spin curves, 
there
exist a representative $M\in\Pic X$ of $\lambda$, a one parameter 
smoothing
$f\colon\X\to S$ of $X$,
and a line bundle $\ma{M}$ on $\X$, such that $\ma{M}\otimes
\ol_{X}\simeq M$
and such that
$$
(\ma{M}^*)^{\otimes r}\simeq (\omega _{f}^*)^{\otimes l}.
$$
Then we have
$$
\ma{M}^{\otimes r} \simeq \omega _{f}^{\otimes l} \otimes\ol_{\X}
(D)
$$
where $D$ is a Cartier divisor supported on the central fiber.
Finally, the line bundle
$T={\mathcal O}_{\X}(D)\otimes {\mathcal O}_{X}$ is a twister of
$X$; therefore, restricting the above relation to $X$, we are done.

Now we prove the converse: let $X$ be a quasistable curve of genus $g$
and let
$T\in\Tw(X)$
be such that
$(\un{\deg}\ T+l\un{w}_X)r^{-1}$ is integer and balanced.
Observe that $\omega_X^{\otimes l}\otimes T$ has degree divisible by
$r$ on each irreducible component of $X$, so it admits some $r$-th 
root
$M\in\Pic X$. Let us show that $M\in\sgtil$.

$M$ is balanced, so it defines a point in $\pdg$.
Moreover, by definition of twister, there exist
a 1-parameter smoothing $f\colon\ma{X}\to S$
of $X$ and a Cartier divisor $D$, supported on the central fiber, such
that
\[ \ol_{\ma{X}}(D)\otimes\ol_X=T.\]
Then $M$ is an $r$-th root of
$(\omega_f^{\otimes l}\otimes\ol_{\X}(D))\otimes\ol_X$, and by
\ref{fibrati}
it has a (unique) extension to $\ma{M}\in\Pic\ma{X}$ such that
$$\ma{M}^{\otimes r}\simeq \omega_f^{\otimes l}\otimes\ol_{\X}(D).$$
Hence $M$ is a limit of smooth
$r$-spin curves, and, as such, it lies in $\sgtil$.
\end{proof}
\begin{remark}
\label{comb}
The proposition
highlights the combinatorial nature of $\sgtil$, by
saying that
the points of $\sgtil$
are determined by all twisters $T$ on $X$ such
that
$r^{-1}\un{\deg}\ (T\otimes \omega _X^{\otimes l})$ is integer and
balanced.
It is easy to see that this is equivalent to the following properties
of $\un{\deg}\ T$:
\begin{enumerate}[{\it (i)}]
\item
\label{exc}
$\deg_{E_i}T=r$ for every exceptional component $E_i$ of $X$;
\item
\label{lw}
$\deg_{X_j}T\equiv (-lw_j)\pmod r$ for every non exceptional
component $X_j$ of $X$;
\item $|\deg_Z T|\leq \dfrac{k_Zr}{2}$ for every subcurve $Z$ of $X$.
\end{enumerate}
\end{remark}
\subsection{Fiberwise description}
Consider again the embedding $\chi$ defined in \ref{new}.
As $\sg$ is dense in $\sgbar$ (and in $\sgtil$ by definition),
$\chi$ will be viewed as
a birational map
$$ \chi\colon\sgbar\dasharrow\sgtil.$$

For any stable curve $C$ of genus $g$, we denote by $\sC$ the fiber of
$\widehat{\ph}\colon\sgtil\to\mgbar$ over $[C]$.
Throughout this section we assume that the automorphism group of $C$ 
is trivial.
If
$C$ is smooth, then $\sC=\sbarC$ is the reduced, zero-dimensional
scheme parametrizing $r$-th roots of
$\omega_C^{\otimes l}$.

We have a first, easy result relating the two above moduli spaces:
\begin{lemma}
\label{ballimit}
Let $(X,L,\alpha)\in\sgbar$ be a limit root and let $C$ be a stable 
curve.
\begin{enumerate}[(i)]
\item
The map $\chi$ is regular at every point $(X,L,\alpha)$ such that
$L$ is balanced.
\item
\label{dim2}
If every limit root in $\sbarC$ is
balanced (for instance, if $C$ is irreducible), then the restriction
$\chi_{|\sbarC}\colon\sbarC
\to\sC$ is surjective. In particular, $\dim \sC=0$.
\end{enumerate}
\end{lemma}
\begin{proof}
The argument for $(i)$ is analogous to the one proving \ref{BI}.
The requirement that $L$ be balanced naturally defines a subfunctor of
$\ov{\mathcal{S}}^{\,r,\,l}_{g}$, which has a natural inclusion
in the functor $\ov{{\mathcal P}}_{l(2g-2)/r,\,g}$.
Hence the statement
follows from Proposition 8.1 of \cite{caporaso}.

Now $(ii)$. We just proved that $\chi$ is
defined on $\sbarC$ and that, for any $\xi\in\sbarC$ with 
$\xi=(X,L,\alpha)$,
the image $\chi(\xi)$ is the class defined by $L$ in $\sC$.
Let $\lambda\in\sC$ and
choose a representative $M\in\Pic X$ for $\lambda$, where
$X$ is some blow-up of $C$; pick a
one-parameter
smoothing of $X$,
$f\colon\X\to S$, and $\ma{M}\in\Pic\X$ such that
$$
(\ma{M}^*)^{\otimes r}\simeq (\omega _{f^*})^{\otimes l}
$$
(``$*$" means away from the special fiber).
Then by Proposition \ref{ssred}, up to replacing $S$ by a finite 
(ramified)
covering, the family
$(f^*\colon\X^*\to S^*,\ma{M}^*)$ determines a unique limit root
$\xi\in \sbarC$. By continuity, $\chi(\xi)=\lambda$, whence
surjectivity. The fact that {\it (\ref{dim2})} applies to irreducible 
curves descends from Lemma \ref{bal}, {\it (\ref{balirr})}.
\end{proof}
Observe that, if $(X,L,\alpha)$ and $(X,L,\beta)$ are two balanced 
limit
roots of
$\omega_C^{\otimes l}$, then they are identified by 
$\chi$ to the
same point in $\sgtil$, regardless of their being 
isomorphic or not as limits
roots (cf. \ref{pull}, part {\it (ii)}).

The main difference between $\sgtil$ and $\sgbar$ is that,
if $r\geq 3$, the morphism
$\widehat{\ph}$ is not finite. First we have
\begin{lemma}
\label{res}
Let $Z$ be an irreducible component of $\sC$.
Then $$\dim
Z\leq b_1(\g_C).$$
In particular, if $C$ is of compact type, then $\dim \sC =0$.
\end{lemma}
\begin{proof}
Let $Z$ be an irreducible component of $\sC$ and let $\lambda \in Z$.
By Proposition~\ref{stil}, there exists a
representative $L\in \Pic X$ for $\lambda$ such that
$L^{\otimes r} \simeq \omega _X
^{\otimes l}\otimes T$ for a twister $T\in \Tw (X)$.
For every given $T$ there are obviously finitely many $L\in \Pic X$
satisfying this relation.
We claim that the twisters $T$ vary in algebraic varieties
of dimension at
most
$b_1(\g _C)$.

Let $\nu\colon X^\nu\to X$ be the normalization, $T$ a twister, and 
set $G:=\nu ^* T\in \Pic \Xn$; then every $T'\in \Tw(X)$ such that
$\nu ^* T'=G$ has the same multidegree as $T$, and hence
all the $r$-th roots of $\omega _X ^{\otimes l}\otimes T'$ lie in
$\sC$ (see Remark \ref{comb}).
Under the map $\nu ^*\colon\Pic X \to \Pic \Xn$
the fiber of $G$ is a $(\C ^*)^b$,
where
$$
b = b_1(\g _X)= b_1(\g _C);
$$
now, $T'$ moves in
such a fiber, and hence in a variety of dimension at most $b_1(\g 
_C)$.
On the other hand, the varying of $G$ does not contribute to the
dimension of $Z$: it is
in fact easy to see that the set
$\nu ^*(\Tw X)$ (where $G$ ranges) is a discrete subset of $\Pic \Xn$.
Therefore we are done.

Finally,
if $C$ is of compact type then $b_1(\g _C)=0$, so the last statement 
is a
special case of the first.
\end{proof}
\begin{prop}
\label{riass}
Let $C=C_1\cup C_2$ with $C_i$ smooth and $\#C_1\cap C_2=k\geq 2$. 
Then
$\dim\sC=k-1$, with the following list of exceptions, where
$\dim\sC=0$:
\begin{enumerate}[{\it (i)}]
\item $r=2$;
\item $k=2$, $l\un{w}_C\equiv 0\pmod r$;
\item $k=3$, $r=4$, $l\un{w}_C\equiv (2,2)\pmod 4$;
\item $k=4$, $r=3$, $l\un{w}_C\equiv 0\pmod 3$.
\end{enumerate}
\end{prop}
\begin{proof}
In the four exceptional cases listed in the statement, it is easy to 
check that all limit
roots in
$\sbarC$ are balanced, hence by Lemma \ref{ballimit}
part~$(\ref{dim2})$ we have $\dim\sC=0$.

Denote by $n_1,\dotsc,n_k$ the nodes of $C$, and by $\nu\colon 
C^\nu\to C$ the normalization.
We divide the proof in three steps.

\

\noindent
{\bf Step 1.} {\it
Fix a twister $T$ on $C$. Let $T'\in\Pic C$ be such that
$\nu^*(T')=\nu^*(T)$; then $T'$ is a twister.}

Let $\ma{C}\to S$ be such that
$T=\ol_{\ma{C}}(\sum_ia_iC_i)\otimes\ol_C$,
and suppose that $\ma{C}$ has an $A_{h_j-1}$ singularity
in $n_j$.
Introduce the universal deformation $\un{\ma{C}}\to\un{D}$ of $C$, 
where
$\un{D}$ is a smooth polydisc. Assume, as usual, that the locus where
the nodes $n_1,\dotsc,n_k$ are preserved is defined by the vanishing 
of the first $k$ coordinates
on $\un{D}$. Consider the
family $\ma{C}'\to S$ obtained by pulling back $\un{\ma{C}}$ via the 
morphism
\begin{align*}
S& \la \un{D}\\
t&\mapsto (c_1t^{h_1},\dotsc,c_k t^{h_k},0,\dotsc,0),
\end{align*}
where $t$ is the local coordinate on $S$ and the $c_j\in \C^*$ are 
arbitrary.
Then $\ma{C}'$, just like $\ma{C}$, has an $A_{h_j-1}$ singularity
at $n_j$; therefore $T' :=\ol_{\ma{C}'}(\sum_ia_iC_i)\otimes \ol_C$ 
is a
line bundle and thus a
twister. Finally,
is easy to see that $T'$
has gluing datum $c_j$ over $n_j$.
The conclusion is that, for any gluing assignment of $\nu ^*T$ over 
the nodes of $C$, we obtain a twister, as desired.

\

From now on we exclude cases $(i), \dotsc , (iv)$ above.

\noindent
{\bf Step 2.} {\it There exists a
nonzero twister $T_0$ on $C$ such that the $r$-th roots of
$\omega _C^{\otimes l}\otimes T_0$ are stably balanced and,
in particular, correspond to points of
$\sC$.}

Pick $k$ non zero integers $b_1,\dotsc,b_k$ having
the same sign.
As usual, denote by $p_1,\dotsc,p_k$
(respectively,
$q_1,\dotsc,q_k$) the inverse images in $C_1$
(respectively, in $C_2$) of the nodes.
Define $G\in \Pic C^{\nu}$ by

$$G_{|C_1}= \ol_{C_1}
(\sum_{i=1}^kb_ip_i),\quad G_{|C_2}=
\ol_{C_2}
(-\sum_{i=1}^kb_iq_i).$$
Then there exists a twister $T$ on $C$ such that $\nu ^* T =G$. In 
fact,
set $h:=\operatorname{lcm}(|b_1|,\dotsc,|b_k|)$ and
$h_i:=h/|b_i|$. Consider a one-parameter smoothing $\ma{C}\to S$ of 
$C$
having an $A_{h_i-1}$ singularity in $n_i$ for $i=1,\dotsc,k$.
Set $D:=h C_2$ if the $b_i$ are positive, $D:=-h C_2$ if the $b_i$
are negative. Then $D$ is Cartier and
$T:=\ol_{\ma{C}}(D)\otimes \ol_C$ is a twister such that
$T_{|C_1}=\ol_{C_1}(b_1p_1+\cdots+b_kp_k)$ and
$T_{|C_2}=\ol_{C_2}(-b_1q_1-\cdots-b_kq_k)$, that is, $\nu ^* T =G$.

Now set $s = \sum b_i$; obviously $|s|\geq k$.
By Proposition \ref{stil}, Step 2 will be proved by solving the 
following numerical
problem.

Find an integer $s$ such that
\begin{enumerate}[$(1)$]
\item
$
|s|\geq k
$
\item $lw_1+s \equiv 0 \pmod r\ $ and \ $lw_2-s \equiv 0 \pmod r$
\item $|s |< \dfrac{kr}{2}$
\end{enumerate}
where (2) and (3) are parts {\it (ii)} and {\it (iii)} of
Remark~\ref{comb}, and we ask for strict inequality in (3) because we 
want
the degree to be stably balanced.

Write $s=rj-lw_1
$ for some integer $j$.
It finally suffices to find an integer $j$ so that (1) and (3) above 
are satisfied, that is
$
k\leq |rj-lw_1|<\frac{k r}{2}
$. In other words, we are reduced to showing that the set
$$
\left(-\frac{k}{2}+\frac{lw_1}{r}\,,\,
-\frac{k}{r}+\frac{lw_1}{r}\right]\cup
\left[\frac{k}{r}+\frac{lw_1}{r}\,,\,\frac{k}{2}+\frac{lw_1}{r}\right)
$$
contains an integer $j$ (provided that cases $(i),\dotsc,(iv)$ are 
excluded).
The proof of this fact is elementary and easy, therefore we omit it.

\

\noindent
{\bf Step 3.} {\it Steps 1 and 2 imply the proposition.}

Pick a twister $T_0$ as in Step 2.
By Proposition \ref{stil}
and Remark \ref{comb}, the set of all such $T_0$ is characterized by 
some
numerical conditions on their multidegree, hence on the multidegree of
their pull-back to the normalization of $C$.
Therefore, by Step 1, if there exists one $T_0$, then there exists a
$b_1(\Gamma
_C)$-dimensional family $\mathcal F$ of them:
$\mathcal F\simeq (\C ^*)^{b_1(\Gamma _C)}$ is the entire fiber of 
$\nu ^*$ over $\nu
^*(T_0)$. This yields that for every
$T\in
\mathcal F$ the $r$-th roots of $\omega _C^{\otimes l}\otimes T$
lie in $\sC$, and correspond to different points, because they are 
stably
balanced.
We conclude $\sC$ contains subvarieties of dimension $b_1(\Gamma
_C)=k-1$; by \ref{res}
we
are done.
\end{proof}
\subsection {Comparison results}
We begin with the
case $r=2$, which is essentially already known, as we will explain.
Under this condition limit roots are always balanced:
\begin{lemma} Assume that $r=2$.
Let $C$ be a stable curve and let $(X,L,\alpha)\in 
\ov{S}^{\,2,\,l}_g$ be a
limit square root of
$(C,\omega_C^{\otimes l})$.
Then $L$
is balanced
and the orbit of the corresponding Hilbert point is closed in the 
GIT-semistable
locus.
\end{lemma}
\begin{proof}
Fix a subcurve $Z$ of $X$, let $k_Z$ and $w_Z$ be as usual, and let
$k'_Z$ be the number of exceptional components of $X$ intersecting $Z$
and not contained in $Z$.
Then we have
$$k'_Z\equiv lk_Z\ \pmod 2,$$ $k'_Z\leq k_Z$ and
$2d_{Z}=lw_Z-k'_Z$. The Basic Inequality (see \ref{defBI}) for $Z$
simplifies to
$k'_Z\leq k_Z$, which is always satisfied, thus $L$ is balanced and 
the Hilbert
point is GIT-semistable
(Theorem~\ref{BI}).
Moreover, if equality
holds, every node on $Z\cap \ov{X\smallsetminus Z}$ is exceptional. 
Then Lemma 6.1 of
\cite{caporaso} applies, giving that the Hilbert point has closed 
orbit, as stated.
\end{proof}
In view of \ref{ballimit}, one obtains that
$\chi\colon\ov{S}^{\,2,\,l}_g\dasharrow\widehat{S}^{\,2,\,l}_g$ is 
regular for all $l$; moreover, $\chi$ is bijective because of the 
last assertion of the previous lemma (cf. \cite{fontanari}).

\bigskip
When $r\geq 3$ the situation is quite different for, in 
general, neither the birational map
$\chi\colon\sgbar\dasharrow\sgtil$ nor its inverse are regular, as we 
shall
presently see.
\begin{thm}
\label{comp} Let $r\geq 3$.
\begin{enumerate}[(i)]
\item
\label{compiso}
$\chi$ is regular in codimension 1.
\item
\label{compchi}
If $r$ does not divide $2l$, then
$\chi$ is not regular in codimension 2.

If $r$ divides $2l$, then
$\chi$ is not regular in codimension 3, unless $r=4$ and $l\equiv
2\pmod 4$.

If $r=4$ and $l\equiv 2\pmod 4$, then $\chi$ is not regular in
codimension 4.
\item\label{chiinv}
If $r\geq 5$ and $g\geq 7$, or $\,r=3,4$ and
$g\geq 10$, then
$\chi ^{-1}$ is not regular in codimension 1.
\end{enumerate}
\end{thm}
\begin{proof}
First of all, recall that $\ov{\ph}\colon \sgbar\to \mgbar$
is finite. By Lemma~\ref{ballimit}, $\chi$ is regular over the locus 
of
irreducible curves in $\mgbar$. Let $U\subset\sgbar$ be the inverse
image of the locus of curves of compact type in $\mgbar$. Then,
by Proposition \ref{smooth}
applied to $B=\mgbar$, $U$ is normal, so $\chi$ is regular in
codimension 1 on $U$, hence regular in codimension 1 on $\sgbar$. This
settles {\it (i)}.

For $(\ref{compchi})$ we use Proposition~\ref{riass},
which describes boundary strata in $\mgbar$ over which 
$\widehat{\ph}\colon
\sgtil \to \mgbar$ is not finite. Since $\ov{\ph}$
is finite, $\chi$, being dominant, cannot be regular over such strata.

In case $r$ does not divide $2l$
we consider the locus of curves $C=C_1\cup C_2$ with $C_i$
smooth, $\#C_1\cap C_2=2$ and $C_1$ of genus $1$ (hence $C_2$ of genus
$g-2$). Then $l\un{w}_C=(2l,2l(g-2))\not\equiv 0\pmod r$, so by
Proposition~\ref{riass} $\,\widehat{\ph}$ is not finite over this 
locus.

If $r$ divides $2l$ and $(r,l)\neq(4,2+4j)$, again by 
Proposition~\ref{riass},
$\,\widehat{\ph}$ is not finite over all codimension 3 loci
whose general curve is a union of 2 smooth components meeting at 3 
points.

Finally, if $r=4$ and $l\equiv 2\pmod 4$, $\widehat{\ph}$ is not 
finite over
all codimension 4 loci
whose general curve is a union of 2 smooth components meeting at 4
points.

$(\ref{chiinv})$
Consider $\xi=(X,L,\alpha)\in\sbarC$ and a one-parameter smoothing
$(f\colon\X\to S,\ma{L})$ of $\xi$.
Then the moduli morphism
$S^*\to\sgtil$ completes to a morphism $S\to\sgtil$. Denote by
$\lambda\in \sC$ the image of $s_0\in S$. Clearly, if $\chi^{-1}$ is 
defined
at $\lambda$, then $\chi^{-1}(\lambda)=\xi$.
We shall
exhibit, for a suitable choice of $C$, $r$ distinct points
$\xi_1,\dotsc,\xi_r\in\sbarC$ and one-parameter deformations
$f_1,\dotsc,f_r$ as above, such that the associated $\lambda_i$ all
lie in the same irreducible component of $\sC$. This
implies that $\chi^{-1}$ is not defined at $\lambda_i$, as the 
$\xi_i$ are
irreducible components of $\sbarC$.

Let $C=C_1\cup C_2$ with $C_i$ smooth, $C_1\simeq\pr{1}$ and
$\#C_1\cap C_2=k$ (hence $C$ has genus $g\geq k-1$ and $w_1=k-2$).

We treat the case $r\geq 5$, $g\geq 7$ and
$l=1$, the remaining cases being similar. Let $k=7$,
fix 5 nodes $n_1,\dotsc,n_5$
of $C$ and consider the weighted subgraph of
$\g_C$ given by these five nodes with weights $u_1=\cdots=u_5=1$ and
$v_1=\cdots=v_5=r-1$. This weighted graph satisfies conditions (C1)
and (C2) of section~\ref{wg} with respect to $\omega_C$.

Let $X$ be the blow-up of $C$ at $n_1,\dotsc,n_5$.
Fix $L_i\in\Pic C_i$ such that
\[L_1^{\otimes
r}\simeq(\omega_X)_{|C_1}\otimes\ol_{C_1}(-\sum_{j=1}^5p_j)
\quad\text{ and
}\quad L_2^{\otimes
r}\simeq(\omega_X)_{|C_2}\otimes\ol_{C_2}(-\sum_{j=1}^5(r-1)q_j).\]
The choices for the gluings of $L_1$ and $L_2$ over $n_6$ and $n_7$ 
yield
$r$ distinct limit roots
$\xi_i=(X,L_i)$ of $(C,\omega_C)$, $i=\{1,\ldots ,r\}$.
Now consider a general one-parameter smoothing
$(f_i\colon \X\to S,
\ma{L}_i)$ of $(X,L_i)$. The surface $\X$ is smooth at
$q_1,\dotsc,q_5,n_6,n_7$ and has $A_{r-2}$ singularities at
$p_1,\dotsc,p_5$; therefore $C_2$ is a Cartier divisor. The line 
bundle
$\ma{M}_i:=\ma{L}_i\otimes
\ol_{\X}(-C_2)$ has degree zero on all exceptional components of $X$, 
hence
it descends to $\ma{M}_i'\in\Pic\ma{C}$, where
$\ma{C}\to S$ is the stable model of $\X\to S$.
Now set $M_i:=
(\ma{M}_i')_{|C}$. Then $M_i$ is stably balanced and determines a 
point
$\lambda_i\in \widehat{S}^{\,r,\,1}_C$.
Observe that $M_i^{\otimes r}\simeq\omega_C\otimes T_i$, where
$T_i$ is a twister on $C$ and
\[ (T_i)_{|C_1}\simeq\ol_{C_1}(-\sum_{j=1}^5p_j),
\qquad
(T_i)_{|C_2}\simeq\ol_{C_2}(\sum_{j=1}^5q_j). \]
Then, as in the proof of Proposition \ref{riass}, we see that there is
a 6-dimensional irreducible component $Z$ of $\widehat{S}^{\,r,\,1}_C$
containing
$\lambda_1,\dotsc,\lambda_r$.
Moreover, as
we vary the families $f_i$, the points $\lambda_i$ cover
an open subset of $Z$. This shows that $Z$ is entirely contained in
the subscheme where $\chi^{-1}$ is not defined.
Since the curve $C$ itself varies in a codimension 7 stratum of 
$\mgbar$,
we obtain that $\chi^{-1}$ is not regular over a codimension 1 
subscheme of
$\widehat{S}^{\,r,\,1}_g$.
\end{proof}
Our comparison results extend easily to torsion-free sheaves 
compactifications.
More precisely, recall that
Pandharipande constructs in \cite{pandha} a
compactification $\ov{U_g(d,1)}$ of $P_{\,d,g}$ over $\mgbar$, by
means of rank 1 torsion-free sheaves, and shows that this
compactification is isomorphic to $\pdg$ under the correspondence
given by Proposition \ref{sheaves}.
Let $\jar$ be the coarse moduli space of Jarvis's stack
${\text{\scshape Root}}_{g}^{1/r}(\omega^{\otimes l})$
of $r$-th root sheaves of $\omega^{\otimes l}$ over $\mgbar$
(see section \ref{spin}).
Set $d=l(2g-2)/r$.
Then we have a commutative diagram:
\[
{}{}\xymatrix{ {\sgbar}\ar[r]^{\sim\ }\ar@{-->}[d]^{\chi} &
{}{} {\jar}\ar@{-->}[d]^{\chi'}\\
{}{}{\pdg}\ar[r]^{\sim} & {\ov{U_g(d,1)}}}
\]
Therefore all results about $\chi$ hold equally for $\chi'$. In
particular, we have:
\begin{cor}
The coarse moduli spaces of the stacks
${\text{\scshape Root}}_{g}^{1/r}(\omega^{\otimes l})$
and $\ov{\mathfrak{S}}^{\,1/r}_g(\omega^{\otimes l})$
(in the notation of \cite{jarvis2})
do not embed in
$\ov{U_g(d,1)}$.
\end{cor}
\begin{proof}
For $\jar$, the statement is an immediate consequence of Theorem
\ref{Jarvis} and Theorem \ref{comp}, $(\ref{compchi})$.
Let ${T}^{1/r}_g(\omega^{\otimes
l})$ be the coarse moduli space of the stack
$\ov{\mathfrak{S}}^{\,1/r}_g(\omega^{\otimes l})$.
Then ${T}^{1/r}_g(\omega^{\otimes
l})$ is the normalization of $\jar$, so it is finite over $\mgbar$.
The statement follows.
\end{proof}

The literature on compactified Picard varieties is quite rich and varied,
so that it would  be worthwhile to study the same issue with respect to other
compactifications, such as those of \cite{simpson}, \cite{esteves}, \cite{jarvis4},
for example. This naturally leads to the problem of 
comparing among each other various constructions of 
the compactified jacobian;  the general situation
is not as clear-cut as it is for the two compactifications (\cite{caporaso} and
\cite{pandha}) used  in our  paper.
 Recent
results in \cite{alexeev} show that, under suitable hypotheses  the spaces of
\cite{caporaso} and  \cite{simpson} coincide; other correlation results
are, at least to us, not known.

We conclude with a detailed example.
\begin{example}
\label{fibra}
Let $C$ be as in Example \ref{es1}. Assume also, for simplicity, 
that  $C_i\not\simeq \mathbb{P}^1 $ for $i=1,2$.
Under this 
assumption, two limit third roots are isomorphic if and only if
the 
underlying line bundles are isomorphic off the exceptional components.
The quasistable curves having stable model $C$ are $X_i$, obtained
blowing up the node
$n_i$ for $i=1,2,3$, and $X_{j,k}$ and $X_{1,2,3}$, described
in Example \ref{es1}.
The balanced multidegrees are:
\begin{enumerate}[$\circ$]
\item
$\un{d}=(d_{C_1},d_{C_2})=(1,-1)$, $(0,0)$ and $(-1,1)$ on $C$;
\item
$\un{d}=(d_{C_1},d_{C_2},d_{E_i})=(0,-1,1)$ and $(-1,0,1)$ on $X_i$;
\item
$\un{d}=(d_{C_1},d_{C_2},d_{E_j},d_{E_k})=(-1,-1,1,1)$ on $X_{j,k}$.
\end{enumerate}
All these degrees are stably balanced;
no multidegree on $X_{1,2,3}$ is balanced.

Let $\ov{P}_{C}^0$ be
the compactified jacobian of $C$.
Recall from \cite{caporaso}, \S 7.3, that
$\ov{P}_{C}^0$ is the union of strata $V^{\un{d}}_I$ parametrized by
the
balanced multidegrees listed above. Each stratum $V^{\un{d}}_I$ is
isomorphic to the variety of line bundles of multidegree $\un{d}$ on 
the
quasistable curve obtained blowing-up the nodes $\{n_i\}_{i\in I}$ of
$C$ (we denote by $V^{\un{d}}$ the strata $V^{\un{d}}_{\emptyset}$).
The three open strata $V^{(1,-1)}$, $V^{(0,0)}$ and $V^{(-1,1)}$ are
$(\C^*)^2$-bundles over
$\Pic^1C_1\times\Pic^{-1}C_2$, $\Pic^0C_1\times\Pic^{0}C_2$
and $\Pic^{-1}C_1\times\Pic^{1}C_2$ respectively;
their
closures are the irreducible components of $\ov{P}_{C}^0$. In
$\ov{V}^{(1,-1)}$ and $\ov{V}^{(-1,1)}$, the fibers are compactified
by
$\pr{2}$; in $\ov{V}^{(0,0)}$ by $\pr{2}$ blown-up in three
points. These three irreducible components intersect along the
codimension 1 strata; see \cite{caporaso}, page 651 for more details.

Recall that $g^{\nu}=g-2$; we are going to show that $\fc$ consists 
of $15\cdot
3^{2g^{\nu}}$ points and
$2\cdot 3^{2g^{\nu}}$ disjoint copies of $\pr{2}$.

Among the points, $9\cdot3^{2g^{\nu}}$ lie
in $V^{(0,0)}$, and are the $L\in\Pic C$ with $L^{3}\simeq\ma{O}_C$;
the remaining
$6\cdot 3^{2g^{\nu}}$ are the $L\in\Pic X_{j,k}$ with
$L_{|\w{X}_{j,k}}^{\otimes 3}\simeq\ma{O}_{\w{X}_{j,k}}
(-2p_j-q_j-p_k-2q_k)$ or
$L_{|\w{X}_{j,k}}^{\otimes
3}\simeq\ma{O}_{\w{X}_{j,k}}(-p_j-2q_j-2p_k-q_k)$;
they lie in $V_{j,k}^{(-1,-1,1,1)}$.
All these points are limit
roots and lie also in $\ov{S}^{\,3,\,0}_C$; notice that the
remaining points in
$\ov{S}^{\,3,\,0}_C$ have unbalanced multidegree.

Each $\pr{2}$ in $\fc$ corresponds to a choice of $L_1\in\Pic C_1$
and $L_2\in\Pic C_2$ such that either $L_1^{\otimes
3}\simeq\ma{O}_{C_1}(p_1+p_2+p_3)$ and $L_2^{\otimes
3}\simeq\ma{O}_{C_2}(-q_1-q_2-q_3)$, or conversely. We describe the
first case, the other being symmetric. So, fix such $L_1$ and
$L_2$. Notice that $(L_1,L_2)\in\Pic^1C_1\times\Pic^{-1}C_2$. Then
the corresponding $\pr{2}$ is the fiber over $(L_1,L_2)$ of
the $\pr{2}$-bundle
$h\colon\ov{V}^{(1,-1)}\to\Pic^1C_1\times\Pic^{-1}C_2$.
This
$\pr{2}$ contains three types of points.
\begin{center}
{}{}\scalebox{0.45}{\includegraphics{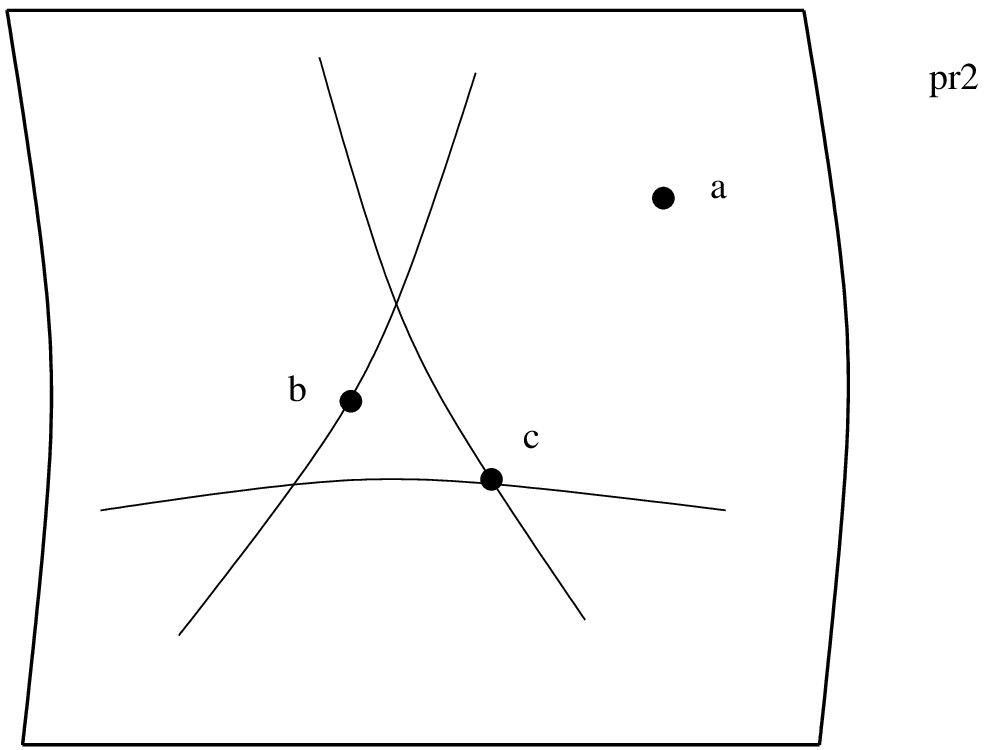}}
\end{center}
The points in
$V^{(1,-1)}$ ($\alpha$ in the figure) correspond to the $L\in\Pic C$ 
with
$L_{|C_1}\simeq L_1$ and $L_{|C_2}\simeq L_2$; the open stratum
$h^{-1}(L_1,L_2)\cap
V^{(1,-1)}=(\C^*)^2\subset\pr{2}$ parametrizes the gluing data for
$L_1$ and $L_2$ over $C$.
These points are obtained
as the central fiber of a 1-parameter family $(\ma{C}\to S,\ma{L})$,
where $\ma{C}$ is smooth and
$\ma{L}\in\Pic\ma{C}$ is
such that $\ma{L}^{\otimes 3}\simeq\ma{O}_{\ma{C}}(-C_1)$.
The line bundle
$\ma{O}_{\ma{C}}(-C_1)_{|C}$ depends on the family.

The points in $V_1^{(0,-1,1)}\cup V_2^{(0,-1,1)}\cup V_3^{(0,-1,1)}$
($\beta$ in the figure)
are the $L\in\Pic X_i$ with $L_{|C_1}\simeq
L_1\otimes\ma{O}_{C_1}(-p_i)$ and
$L_{|C_2}\simeq L_2$.
For each $i$, the open stratum $h^{-1}(L_1,L_2)\cap 
V_i^{(0,-1,1)}\simeq
\C^*\subset\pr{2}$ parametrizes the gluing data for
$L_1\otimes\ma{O}_{C_1}(-p_i)$ and $L_2$ over $\w{X}_i$.
These points are
obtained as the central fiber of a
1-parameter family $(\ma{X}_i\to S,\ma{L})$, where $\ma{X}_i$ has
an $A_2$ singularity in
$q_i$ and $\ma{L}\in\Pic\ma{X}_i$ is
such that $\ma{L}^{\otimes 3}\simeq\ma{O}_{\ma{X}_i}(-C_1-3E_i)$.

Finally, the three points in $V_{1,2}^{(-1,-1,1,1)}$,
$V_{1,3}^{(-1,-1,1,1)}$ and $V_{2,3}^{(-1,-1,1,1)}$ ($\gamma$ in
the figure)
are the $L\in\Pic X_{j,k}$ with $L_{|C_1}\simeq
L_1\otimes\ma{O}_{C_1}(-p_j-p_k)$ and
$L_{|C_2}\simeq L_2$. They are obtained as
the central fiber of a
1-parameter family $(\ma{X}_{j,k}\to S,\ma{L})$,
where $\ma{X}_{j,k}$ has $A_2$ singularities in
$q_j$ and $q_k$, and $\ma{L}\in\Pic\ma{X}_{j,k}$ is
such that $\ma{L}^{\otimes 
3}\simeq\ma{O}_{\ma{X}_{j,k}}(-C_1-3E_j-3E_k)$.
\end{example}

\bigskip

\noindent L.\ Caporaso\\
Dipartimento di Matematica, Universit\`a Roma Tre\\
Largo S.\ L.\ Murialdo, 1\\
00146 Roma - Italy\\
caporaso@mat.uniroma3.it

\bigskip

\noindent C.\ Casagrande\\
Dipartimento di Matematica ``L.\ Tonelli'', Universit\`a di Pisa\\
Largo B.\ Pontecorvo, 5\\
56127 Pisa - Italy\\
casagrande@dm.unipi.it

\bigskip

\noindent M.\ Cornalba\\
Dipartimento di Matematica ``F.\ Casorati'', Universit\`a di Pavia\\
Via Ferrata, 1\\
27100 Pavia - Italy\\
maurizio.cornalba@unipv.it

\end{document}